\documentclass[aop,MSNbibl,seceqn,citesort,dvips]{arximspdf}
\usepackage{graphicx}

\doi{10.1214/11-AOP659}
\volume{39}
\issue{5}
\pubyear{2011}
\firstpage{1768}
\lastpage{1814}

\makeatletter

\newcommand{\mod}{\operatorname{mod}}

\newtheorem{theorem}{Theorem}[section]
\newtheorem{corollary}[theorem]{Corollary}
\newtheorem{lemma}[theorem]{Lemma}
\newproclaim{assumptions}[theorem]{Assumptions}
\newproclaim{question}[theorem]{Question}
\newproclaim{remark}[theorem]{Remark}

\newtheorem{proposition}[theorem]{Proposition}
\newtheorem{conjecture}[theorem]{Conjecture}

\newcommand{\R}{\mathbb{R}}
\newcommand{\C}{\mathbb{C}}
\newcommand{\Z}{\mathbb{Z}}
\newcommand{\diam}{\operatorname{diam}}
\newcommand{\const}{\operatorname{const}}
\newcommand{\area}{\operatorname{area}}
\newcommand{\dist}{\operatorname{dist}}
\newcommand{\eps}{\epsilon}
\newcommand{\E}{\mathbf{E}}

\newcommand{\closure}{\overline}
\newcommand{\Quad}{\mathcal Q}
\newcommand{\Scal}{\mathcal H}
\newcommand{\clo}{\mathfrak F}
\newcommand{\topo}{\mathcal T}
\newcommand{\Pdisc}{\mu}
\newcommand{\Plim}{{\mu_0}}
\newcommand{\dsamp}{\omega}
\newcommand{\bay}{B}
\newcommand{\cross}{{\boxminus}_}

\newcommand{\rswa}{a_1}
\newcommand{\rswpa}{a_4}
\newcommand{\rsww}{\Delta}

\newcommand{\prsw}{P}
\newcommand{\exep}{{S}}
\newcommand{\interface}{\gamma}
\newcommand{\tsum}{ \sum}

\makeatother

\begin{document}
\begin{frontmatter}

\title{On the scaling limits of planar percolation}
\runtitle{On the scaling limits of planar percolation}

\begin{aug}
\author{\fnms{Oded} \snm{Schramm}\thanksref{t1}} and
\author[A]{\fnms{Stanislav} \snm{Smirnov}\corref{}\thanksref{t2}\ead[label=e2]{stanislav.smirnov@unige.ch}%
\ead[label=u2,url]{http://www.unige.ch/\textasciitilde smirnov}}\\[5pt]
(with an Appendix by
\author[B]{\fnms{Christophe} \snm{Garban}\ead[label=e3]{christophe.garban@umpa.ens-lyon.fr}%
\ead[label=u3,url]{http://www.umpa.ens-lyon.fr/\textasciitilde cgarban/}})
\runauthor{O. Schramm and S. Smirnov}
\affiliation{Microsoft Research, Universit\'{e} de Gen\`{e}ve and\break
CNRS, ENS Lyon}
\address[A]{S. Smirnov\\
Section de Math\'{e}matiques\\
Universit\'{e} de Gen\`{e}ve\\
2-4 rue du Li\`{e}vre, CP 64\\
1211 Gen\`{e}ve 4\\
Switzerland\\
\printead{e2}\\
\printead{u2}}
\address[B]{C. Garban\\
CNRS, UMPA\\
ENS de Lyon\\
46, all\'{e}e d'Italie\\
69364 Lyon Cedex 07\\
France\\
\printead{e3}\\
\printead{u3}}
\end{aug}

\thankstext{t1}{December 10, 1961--September 1, 2008.}

\thankstext{t2}{Supported by the European Research Council AG CONFRA,
the Swiss National Science Foundation and by the Chebyshev Laboratory
(Faculty of Mathematics and Mechanics, St. Petersburg State University)
under the grant of the government of the Russian Federation.}

\received{\smonth{2} \syear{2010}}
\revised{\smonth{2} \syear{2011}}

%
\begin{abstract}
We prove Tsirelson's conjecture that any scaling limit of
the critical planar percolation is a black noise.
Our theorems apply to a number of percolation models,
including site percolation on the triangular grid
and any subsequential scaling limit of bond percolation
on the square grid.
We also suggest a natural construction for
the scaling limit of planar percolation,
and more generally of
any discrete planar model
describing connectivity properties.
\end{abstract}

%

\end{frontmatter}

\section{Introduction}

\subsection{Motivation}
This paper has a two-fold motivation:
to propose a~new construction for the (subsequential) scaling limits of the
critical and near-critical percolation in the plane,
and to show that such limits are two-dimensional \textit{black noises}
as suggested by Boris Tsirelson.

Percolation is perhaps the simplest statistical physics model
exhibiting a phase transition.
We will be interested in planar percolation models,
the archetypical examples being
site percolation on the triangular lattice
and bond percolation on the square lattice.
More generally, our theorems apply to a wide class of
models on planar graphs,
satisfying the assumptions
outlined below.

Consider a planar graph and fix a \textit{percolation probability} $p\in[0,1]$.
In \textit{site percolation}, each site is declared \textit{open} or
\textit{closed}
with probabilities $p$ and $1-p$, independently of others.
In the pictures, we represent open and closed sites
by blue and yellow colors correspondingly.
One then studies open \textit{clusters}, which
are maximal connected subgraphs of open sites.
In \textit{bond percolation}, each bond is similarly declared open or
closed independently of others,
and one studies clusters of bonds.

Percolation as a mathematical model was originally introduced by
Broadbent and Hammersley
to model the flow of liquid through a random porous medium.
It is assumed that the liquid can move freely between the
neighboring open sites,
and so it can flow between two regions if there is an open path between them.
Fixing a quad $Q$ ($=$ a topological quadrilateral),
we say that there is an \textit{open crossing}
(between two distinguished opposite sides)
for a given percolation configuration,
if for its restriction to the quad,
there is an open cluster intersecting both opposite sides.
One then defines the crossing probability for $Q$ as the probability
that such a crossing exists.

We will be interested in the (subsequential) scaling limits of
percolation, as the lattice mesh 
tends to zero.
It is expected that for a large class of models the \textit{sharp
threshold} phenomenon
occurs: there is a critical value $p_c$ such that for a fixed quad $Q$
the crossing probability tends to $1$ for $p>p_c$ and to $0$ for $p<p_c$
in the scaling limit.
For bond percolation on the square lattice and site percolation on the
triangular lattice,
it is known that $p_c=1/2$ by the classical theorem of Kesten~\cite{kesten-12}.

In the critical case $p=p_c$ the crossing probabilities are expected
to have nontrivial scaling limits.
Moreover, those are believed to be universal, conformally invariant,
and given by Cardy's formula \cite{langlands-bams,Cardy-92}.
For the models we are concerned with,
this was established only for site percolation on the triangular
lattice \mbox{\cite{smirnov-cras,smirnov-perc}}.
However, the Russo--Seymour--Welsh theory provides
existence of nontrivial subsequential limits,
which is enough for our purposes.
Our theory also applies to near-critical scaling limits like in
\cite{nolin-near},
when the mesh tends to zero and the percolation probability to its
critical value
in a coordinated way, so as the
crossing probabilities have nontrivial limits,
different from the critical one.
The Russo--Seymour--Welsh theory still applies below any given scale.

Having a (subsequential) scaling limit of the (near) critical percolation,
we can ask how to describe it.
The discrete models have finite $\sigma$-fields,
but in the scaling limit the cardinality blows up,
so we have to specify how we pass to a limit.
First, we remark that we restrict ourselves to the \textit{crossing events}
(which form the original physical motivation for the percolation model),
disregarding events of the sort
``number of open sites in a given region is even''
or ``number of open sites in a given region is bigger than the number
of closed ones.''
While such events could be studied in another context,
they are of little interest in the framework motivated by physics.
Furthermore, we restrict ourselves to the \textit{macroscopic}
crossing events, that is, crossings of quads of fixed size.
Otherwise, one could add to the $\sigma$-field crossing events for
microscopic quads (of the lattice mesh scale or some intermediate scale),
which have infinitesimal size in the scaling limit.
The resulting construction would be an extension of ours,
similar to a nonstandard extension of the real numbers.
But because of the locality of the percolation,
different scales are independent, so it won't yield new nontrivial information
(as our results in fact show).
However, we would like to remark
that in dependent models such extensions
could provide new information.
For example, Kenyon has calculated probabilities
of local configurations appearing in domino tilings and loop erased
random walks \cite{kenyon-lectures},
while such information is lost in describing the latter by an SLE curve.

Summing it up, we restrict ourselves only
to the scaling limit of the
\textit{macroscopic crossing events} for percolation---it
describes the principal physical properties
(though adding microscopic events can be of interest).
The question then arises, how to describe
such a scaling limit, and several approaches were proposed.
We list some of them along with short remarks:
\begin{longlist}[(2)]
\item[(1)] \textit{Random coloring of a plane} into yellow and blue colors
would be~a~lo\-gical candidate,
since in the discrete setting we deal with random colorings of graphs.
However, it is difficult to axiomatize the allowed class
of ``percolation-induced colorings,''
and connectivity properties are hard to keep track of.
In particular, in the scaling limit almost every point does not belong
to a cluster,
but is rather surrounded by an infinite number of nested open and
closed clusters,
and so has no well-defined color.
\item[(2)] \textit{Collection of open clusters} as random compact subsets of
the plane would satisfy
obvious precompactness (for Borel probability measures on collections
of compacts endowed with
some version of the Hausdorff distance), guaranteeing existence
of the subsequential scaling limits.
However, crossing probabilities cannot be extracted:
as a toy model imagine two circles, $O$ and $C$, intersecting at two
points, $x$ and $y$.
Consider two different configurations:
in the first one, $O\setminus\{x\}$ is an open crossing and
$C\setminus\{y\}$ is a~closed one;
in the second one, $O\setminus\{y\}$ is an open crossing and
$C\setminus\{x\}$ is a~closed one.
In both configurations, the cluster structure is the same
(with~$O$ being the open and $C$ being the closed cluster),
while some crossing events differ.
Thus one has to add the state of the ``pivotal'' points
(like $x$ and $y$) to the description,
making it more complicated.
Constructing such a limit, proving its uniqueness,
and studying it would require considerable technical work.
\item[(3)] \textit{Aizenman's web} represents percolation configuration by a
collection
of all curves ($=$ crossings) running inside open clusters.
Here pre-compactness follows from the Aizenman--Burchard work
\cite{aizenman-burchard},
based on the Russo--Seymour--Welsh theory, and crossing events are easy
to extract.
Moreover, it turns out that the curves involved are almost surely
H\"{o}lder-continuous, so this provides a nice geometric description of
the connectivity structure.
However, there is a lot of redundancy in taking \textit{all} the curves
inside clusters,
and there are deficiencies similar to those of the previous approach.
\item[(4)] \textit{Loop ensemble}---instead of clusters themselves, it is enough
to look at the interfaces between open and closed clusters, which
form a collection of nested loops.
Such a limit was constructed by Camia and Newman in \cite{camia-newman-cmp}
from Schramm's SLE loops.
However, it was not shown that such a limit is full (i.e., contains all
the crossing events),
and in general, extracting geometric information is far from
being straightforward.
\item[(5)] \textit{Branching exploration tree}---there is a canonical way
to explore all the discrete interfaces:
we follow one, and when it makes a loop, branch into two created components.
The continuum analogue was studied
(in relation with the loop ensembles)
by Sheffield in \cite{sheffield-cle},
and this is perhaps the most natural approach for general random
cluster models
with dependence; see \cite{kemppainen-smirnov-fk3}.
Construction is rather technical,
however the scaling limit can be described by branching SLE curves and so
is well suited for calculations.
\item[(6)] \textit{Height function}---in the discrete setting one can
randomly orient interfaces, constructing
an integer-valued height function, which changes by $\pm1$ whenever an
interface is crossed.
This stores all the connectivity information
(at least in our setup of interfaces along the edges of trivalent graphs),
but adds additional randomness of the height change.
It is expected that after appropriate manipulations
(e.g., coarse-graining and compactifying---i.e.,
projecting $\mod1$ to a unit interval)
the height function would converge
in the scaling limit to a random distribution,
which is a version of the free field.
This is the cornerstone of the ``Coulomb gas'' method of Nienhuis
\cite{nienhuis-1984jsp},
which led to many a prediction.
It is not immediate how to connect this approach to the more geometric ones,
retrieving the interfaces ($=$~level curves) from the random distribution,
but SLE theory suggests some possibilities.
Also it is not clear whether we can implement this approach so that
the extra randomness disappears in the limit,
while the connectivity properties remain.
\item[(7)] \textit{Correlation functions}---similarly to the usual
Conformal Field Theory approach,
instead of dealing with a random field as a random distribution, we can
restrict our attention
to its $n$-point correlations. For example, that would be the usual
physics framework for fields
arising in approach (6), but to implement it mathematically
we would need to reconstruct height lines of a field, based on its
correlations, which is perhaps more difficult than developing
(6) by itself,
and possibly requires some extra assumptions on the field.
A more geometric approach would be to study
for every finite collection of points $\{z_j^k\}_{j,k}$
the probability $F_r(z_j^k)_{j,k}$
that for every $k$ all the balls $\{B(z_j^k,r)\}_j$ are connected
by an open cluster.
To obtain the correlation function, one
takes the double limit,
passing to the scaling limit first and sending $r\to0$ afterwards.
Since the probability of a radius $r$ ball being touched by a
macroscopic cluster
(conjecturally) decays as $r^{5/48}$, one has to normalize
accordingly, considering $F$
as a $\prod_{j,k}{|d z_j^k|}^{5/48}$-form.
The geometric approach seems more feasible,
but along the way one would have to further develop
a version of CFT corresponding to such
``connectivity'' correlations.
\item[(8)] \textit{Collection of quads crossed}---in the discrete setting
one can describe a~percolation configuration by listing
the finite collection of all discrete quads (i.e., topological quadrilaterals)
having an open crossing.
We propose a continuous analogue,
described in detail below.
The percolation configuration space thus becomes
a space of quad collections,
satisfying certain additional properties.
\item[(9)] \textit{Continuous product of probability spaces}, or \textit{noise}---this
approach was advocated by Tsirelson
\cite{tsirelson-nonclassical,tsirelson-lnm,tsirelson-experiment}.
In the discrete setting, the site percolation (on a set $V$ of vertices)
is given by the \textit{product probability space}
$(\Omega_V=\prod_{V}\{\mbox{open},\mbox{closed}\}, \mathcal{F}_V,
\Pdisc_V)$.
Note that when $V$ is decomposed into two disjoint subsets,
we obviously have
%
%
\begin{equation}\label{eqprod}
\mathcal{F}_{V} =
\mathcal{F}_{V_1}\times
\mathcal{F}_{V_2} .
\end{equation}
Straightforward generalization of the product space
to the continuous case does not work:
the product
$\sigma$-field for the space $\prod_{\C}\{\mbox{open},\mbox
{closed}\}$
won't contain crossing events.
Instead Tsirelson proposed to send the percolation measure
space to a scaling limit,
described by a \textit{noise}
or a \textit{homogeneous continuous product of probability spaces}.
Namely, for every smooth domain~$D$ the percolation scaling limit
inside it
would be given by a probability space
$(\Omega_D,\mathcal F_D,\Pdisc_D)$,
which is translation invariant, continuous in $D$,
and satisfies the property (\ref{eqprod}) for a smooth domain $V$
cut into two smooth domains~$V_1$, $V_2$ by a curve.
There is a well-developed theory of such spaces,
but establishing the latter property was problematic
\cite{tsirelson-experiment}.
\end{longlist}
To most constructions, one
can also add the information about the closed clusters,
though in principle it should be retrievable from the open ones
by duality.
This list is not exhaustive, but it shows the variety of
possible descriptions of the same object.
Indeed, in the discrete settings the approaches above
contain all the information about macroscopic open clusters,
and are expected to keep it in the scaling limit;
on the other hand with an appropriate setup
[one has to be careful, especially with (6)],
we do not expect to pick up any
extra information on the way.
However, showing that most of these approaches
lead to equivalent results,
and in particular to isomorphic
$\sigma$-fields, \textit{is far from easy}.

We decided to proceed with the approach (8), since it
follows the original physical motivation
for the percolation model,
and also provides us with enough pre-compactness to
establish existence of subsequential scaling limits
\textit{before} we apply any percolation techniques.

It also serves well our second purpose:
to provide escription of the percolation (subsequential) scaling limits,
following Tsirelson's approach (9),
and show that it is a \textit{noise}, that is,
a homogeneous continuous (parameterized by an appropriate Boolean algebra
of piecewise-smooth planar domains)
product of probability spaces.

The question whether critical percolation scaling limit is a noise
was posed by Tsirelson in
\cite{tsirelson-nonclassical,tsirelson-lnm,tsirelson-experiment},
where he has also noted that such noise must be nonclassical, or \textit{black}
(compared to the classical white noise).
The theory of black noises was started by Tsirelson and Vershik in
\cite{tsirelson-vershik}.
In particular, they constructed the first examples of zero and
one-dimensional black noises
(and more of those were found since);
we provide the first genuinely two-dimensional example.

To establish that percolation leads to noise,
it is enough to prove that its scaling limit,
restricted to two adjacent rectangles,
determines the scaling limit in their union.
Following is a version of this result,
perhaps more tangible
(though stated slightly informally here).
It is proved in Proposition \ref{pmain},
which is the most technical result of our paper.
\textit{Consider a rectangle},
\textit{and a smooth path $\alpha$ cutting it.
We show that for every $\eps>0$ there is a finite number of
percolation crossing events in quads disjoint from $\alpha$},
\textit{from which one can reliably predict whether or not the rectangle is crossed,
in the sense that the probability for a mistake is less than $\eps$.}
The essential point is that the set of crossing events that one looks at
\textit{does not depend on the mesh of the lattice}, which allows us to
deduce the needed result,
thus showing that any subsequential scaling limit of the critical percolation
is indeed a \textit{black noise}.

\subsection{Percolation basics and notation}
\label{condperc}

For completeness of exposition, we present most of the needed
background below;
interested readers can also consult the books
\cite{kesten-book,grimmett-book,bollobas-riordan-book,werner-pcmi}.

We start with a rather general setup, which in particular includes all
planar site and bond models
(as well as percolation on planar hypergraphs; cf.
\cite{bollobas-riordan-selfdual}).
Since we will be interested in scaling limits, it is helpful
from the very beginning to work with subsets of the plane, rather than graphs.

We consider a locally finite \textit{tiling} $H$ of the plane (or its
subdomain) by
topologically closed polygonal \textit{tiles} $P$.
Furthermore, we ask it to be trivalent,
meaning that at most three tiles meet at every vertex
and so any two tiles are either disjoint or share a few edges.
A \textit{percolation model} fixes a percolation probability $p(P)\in
[0,1]$ for every tile $P$,
and declares it open (or closed) with probability $p(P)$
[or $1-p(P)$] independently of the others.
In principle, we can work with a random tiling $H$,
with percolation probabilities $p(P,H)$
satisfying appropriate measurability conditions
(so that crossing events are measurable).

More generally, we can consider a \textit{random coloring}
of tiles into open and closed ones, given
by some measure $\Pdisc$ on the
space $\prod_{P\in H}\{\mbox{open},\mbox{closed}\}$
with a product $\sigma$-field
(which contains all events concerned with a finite number of tiles).
Percolation models correspond
to product measures $\Pdisc$.

We define open (closed) \textit{clusters} as
connected topologically closed subsets of the plane---components
of connectivity of the union of open (closed) tiles taken
with boundary.

%
\begin{figure}

\includegraphics{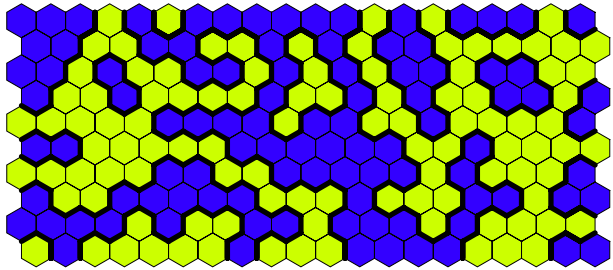}

\caption{Critical site percolation on the triangular lattice. Each
site is represented by a~hexagonal
tile, which is open (blue) or closed (yellow) with probabilities $1/2$
independently of others.
Interfaces between open and closed clusters are pictured in bold.}
\label{fighex}
\end{figure}

All site and bond percolation models on planar graphs
can be represented in this way.
Indeed, given a planar site percolation model,
we replace sites by tiles,
so that two of them share a side if and only if two corresponding
sites are connected by a bond.
It is easy to see that such a collection always exists
(e.g., tiles being the faces of the dual graph;
see Figure \ref{fighex} for the triangular lattice case).
Each tile is declared open (or closed) with the same~probabili\-ty~$p$
(or $1-p$) as the corresponding site.

If our original graph is a triangulation,
then the dual one is trivalent and only three tiles can meet at a point.
However, when our graph has a~face with four or more sides,
two nonadjacent sites can correspond to tiles sharing a~vertex (but not
an edge).
To resolve the resulting ambiguity we insert a~small tile $P$ around
every tiling vertex, where four or more tiles meet,
and set it to be \textit{always closed}, taking $p(P)=0$.
After the modification, at most three tiles can meet at a point,
as required.

\begin{figure}

\includegraphics{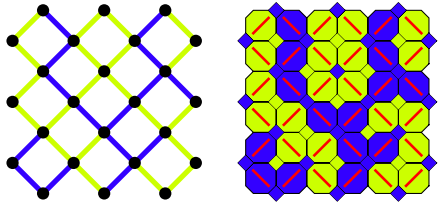}

\caption{Left: critical bond percolation on the square lattice,
every bond is open (blue) or closed (yellow) with probabilities $1/2$
independently of others.
Right: representation by critical percolation on the bathroom tiling,
with direction of original bonds marked on the corresponding octagonal tiles.
Rhombic tiles alternate in color,
with blue (always open) corresponding to
the sites of the original square lattice,
while yellow (always closed)---to the faces.
Octagonal tiles correspond to the bonds, and are open or closed
independently, with probabilities $1/2$.}
\label{figsquare}
\end{figure}

Bond percolation can also be represented by partitioning the plane into tiles.
There will be a tile for every bond, site, and face of the original
bond graph,
chosen so that two different face-tiles or two different site-tiles are
disjoint.
Moreover, a bond-tile shares edges with two site-tiles and two face-tiles,
corresponding to adjacent sites and faces.
The bond-tiles are open (or closed) independently, with
the same probability $p$ (or $1-p$) as the corresponding bonds.
The site-tiles $P$ are then always open [$p(P)=1$], while the
face-tiles $P$ are always closed [$p(P)=0$].
With such construction, no ambiguities arise, and at most three tiles
can meet at a point.
A possible construction for the bond percolation
(with octagons being the bond-tiles) is illustrated in Figure \ref{figsquare}.

%


Following the original motivation of Broadbent and Hammersley,
we represent holes in some random porous medium by open tiles.
A liquid is assumed to move freely through the holes,
thus it can flow between two regions if there is a path of open tiles
between them.
More precisely, we say that there is an \textit{open} (\textit{closed})
\textit{crossing}
between two sets
$K_1$ and $K_2$ inside the set $U$,
if for the percolation configuration restricted to $U$, there is an
open (closed)
cluster intersecting both $K_1$ and $K_2$.
We will be mostly interested in crossings of quads ($=$ topological
quadrilaterals).

Note the following duality property:
\textit{there is an open crossing between two opposite sides of a quad
if and only if
there is no dual} (\textit{i.e.}, \textit{closed}) \textit{crossing
between two other sides}.
Exchanging the notions of open and closed, we conclude
that, on a given tiling, percolation models with probabilities $p(P)$
and $p^*(P)=1-p(P)$ are dual to each other.
In particular, the model with $p(P)\equiv1/2$ is self-dual,
which can often be used to show that it is critical, like for the
hexagonal lattice.
The bathroom-tile model of the square lattice bond percolation
in Figure \ref{figsquare} is dual to itself translated by one lattice step
(so that always open rhombi are shifted to always closed ones),
which also leads to its criticality, \cite{kesten-12}.
There are, however, dualities
(e.g., arising from the star-triangle transformations)
with more complicated relations between~$p$ and $p^*$,
and so there are critical models with $p$ away from $1/2$;
see, for example, \cite{wierman-hex,bollobas-riordan-selfdual}.

Though we have in mind critical models,
we do not use criticality in the proof,
rather some crossing estimates,
which are a part of the Russo--Seymour--Welsh theory.
Those also hold (below any given scale)
for the near-critical models, like those discussed in \cite{nolin-near}.
In a sense the estimates we need should hold
whenever the crossing probabilities do not tend to a trivial limit;
see Remark \ref{remrswnontriv}.

We work with a collection of percolation models $\Pdisc_\eta$
(and sometimes more general random colorings)
on tilings $H_\eta$,
indexed by a set $\{\eta\}$ (e.g., square lattices of different mesh),
and denote by $\Pdisc_\eta$ also the corresponding measure on
percolation configurations
($=$ random coloring of tiles).
By $|\eta|\in\R_+$, we will denote a ``scale parameter'' of the
model $\Pdisc_\eta$.
We assume that all tiles have diameter at most $|\eta|$
(for a percolation on random tiling,
our proof would work under the relaxed assumption
that for any positive $r$, the probability to find in a bounded region
a tile
of diameter bigger than~$r$ tends to zero as $|\eta|\to0$).
This parameter also appears in the Russo--Seymour--Welsh-type estimates
below, which have
to hold on scales larger than $|\eta|$.
Note that in most models $|\eta|$ can be taken to be the supremum
of tile diameters
(or the lattice mesh for the percolation on a lattice).

We will be interested in scaling limits as the scale parameter tends to zero,
looking at sequences of percolation models with $|\eta|\to0$,
which do not degenerate in the limit.
We will show in Remark \ref{remseqnet},
that in our setting working with sequences and nets yields the same results.
As discussed above, all critical percolation models
are conjectured to have a universal and conformally invariant scaling limit.
Moreover, its SLE description allows to calculate exactly many scaling
exponents and dimensions
and we will use this conjectured picture in informal discussions below.
However, to make our paper applicable to a wider class of models,
we will work under much weaker assumptions.

The famous Russo--Seymour--Welsh theory provides universal bounds
for crossing probabilities of the same shape, superimposed
over lattices with different mesh.
Originally established for the critical bond percolation on the square lattice,
it has been since generalized to many critical
(and near-critical, when the percolation probabilities
tend to their critical values at appropriate speeds
as the mesh tends to zero)
models with several different approaches; see
\cite{kesten-book,grimmett-book,bollobas-riordan-book,werner-pcmi}.
We will use one of the recurrent techniques: conditioning on the lowest
crossing.
If a quad is crossed horizontally, one can select the lowest possible
crossing (the closest to the bottom side),
which will then depend only on the configuration below.
If we condition on the lowest crossing, the configuration above it is unbiased,
allowing for easy estimates of events there.
Besides using this neat trick,
we will assume certain bounds on the crossing probabilities for the annuli.

A prominent role is played by
the so-called \textit{$k$-arm} crossing events in annuli,
since they control dimensions of several important sets.
In particular,
for a given percolation model $\Pdisc_\eta$ we denote by $\prsw
^{1}_{\eta}(z,r,R)$
the probability of a \textit{one-arm} event, that
\textit{there is an open crossing connecting two opposite circles of the
annulus $A(z,r,R)$}.
Then this is roughly the probability that the disc $B(z,r)$ intersects
an open cluster of size $\approx R$,
which is expected for $|\eta|<r$ to have a power law like
%
\[
\prsw^{1}_\eta(z,r,R) = (r/R)^{\rswa+o(1)} \qquad\mbox{as } r/R\to
0 ,
\]
morally meaning
(modulo some correlation estimates)
that percolation clusters in the scaling limit
have dimension $2-\rswa$.

Similarly, denote by $\prsw^{4}_\eta(z,r,R)$
the probability of a \textit{four-arm} event, that
\textit{there are alternating open--closed--open--closed
crossings connecting two opposite circles of the annulus $A(z,r,R)$}.
Then this is roughly the probability that changing
the percolation configuration on the disc $B(z,r)$
changes the crossing events on the scale $\approx R$.
Indeed, if all the tiles in $B(z,r)$ are made open,
that connects two open arms; while
making them closed connects two closed arms, destroying the open connection.
For $|\eta|<r$ a power law
%
\[
\prsw^{4}_\eta(z,r,R) = (r/R)^{\rswpa+o(1)} \qquad \mbox{as } r/R\to
0
\]
is expected, roughly meaning that the
\textit{pivotal} tiles (such that altering their state
changes crossing events on large scale) have
dimension $2-\rswpa$ in the scaling limit.

Whenever the RSW theory applies, the probabilities for $k$ arm events
have power law bounds from above and below as $r\to0$.
In fact, the probabilities are conjectured to satisfy
universal (independent of a particular percolation model)
power laws with rational
powers, predicted by the Conformal Field Theory;
see discussion in \cite{smirnov-werner}.
Those predictions have been so far proved for
the critical site percolation on triangular lattice only
\cite{smirnov-werner,lsw-onearm},
giving $\rswa=5/48$ and $\rswpa=5/4$.

Other $k$-arm probabilities are also of interest,
but estimating just the mentioned two allows to apply our methods.
Roughly speaking, we need to know that the percolation
clusters have dimension smaller than $2$
and that pivotal points
have dimension smaller than $1$,
or, in other words, that $\rswa>0$ and $\rswpa>1$.
Note that we do not actually need the power laws,
but rather weaker versions of the corresponding upper estimates.

Namely, we need the following two estimates
to hold \textit{uniformly} for
percolation models under consideration.
\begin{assumptions}
\label{assursw}
There exist positive functions $\rsww_1(r,R)$ and $\rsww_4(r,R)$,
such that
\[
\lim_{r\to0}\rsww_j(r,R)=0 \qquad\mbox{for any fixed } R<R_0 ,
\]
and the following estimates hold
whenever $z\in\C$ and $0<|\eta|<r<R<R_0$.
The probability of one open arm event and
the probability of a similar\vadjust{\goodbreak} event for one closed arm
satisfy
%
%
\begin{equation}
\label{eqwrsw}
\prsw^{1}_{\eta}(z,r,R) \le\rsww_1(r,R) ,
\end{equation}
while the probability of a four arm event satisfies
%
%
\begin{equation}
\label{eqpivo}
\prsw^{4}_{\eta}(z,r,R) \le\biggl(\frac{r}{R}\biggr) \cdot\rsww_4(r,R).
\end{equation}
\end{assumptions}

For $r\ge R$, when the annulus is empty, we set
$\rsww_j(r,R):=1$, so that the function is defined for all positive arguments.
Without loss of generality, we can assume that functions $\rsww_j(r,R)$
are increasing in $r$ and decreasing in $R$.
\begin{remark}[(A stronger RSW estimate)]
\label{remsrsw}
For most percolation models
where the RSW techniques work,
their application would
prove a stronger estimate, with $\rsww_1(r,R)$
replaced in (\ref{eqwrsw}) by $C (r/R)^{\rswa}$ with $\rswa>0$.
However, we can imagine situations where
only our weaker version can a priori be established.
\end{remark}
\begin{remark}[(RSW estimates and scaling limits)]
\label{remrswnontriv}
For most percolation models where the
RSW techniques work,
the estimate (\ref{eqwrsw})
would be equivalent to saying that crossing probability for
any given quad (or for all quads) is uniformly bounded away from $0$
and $1$.
So morally our noise characterization would apply whenever
one can speak of nontrivial subsequential scaling limits
of crossing probabilities.
\end{remark}

The first assumption (\ref{eqwrsw})
is used several times throughout the paper.
In a stronger form, it was one of the original Russo--Seymour--Welsh results,
and has since been shown to hold for a wide range of percolation
models; see
\cite
{bollobas-riordan-book,bollobas-riordan-selfdual,grimmett-book,kesten-book}.
For site percolation on the triangular lattice,
one can take $\rsww_1(r,R)\approx(r/R)^{5/48}$; see \cite{lsw-onearm}.

The second assumption (\ref{eqpivo})
is used only once, namely in the proof of Theorem~\ref{tdiscrete}.
It is known [in a sharp form with $\rsww_4(r,R)\approx(r/R)^{1/4}$]
for site percolation on the triangular lattice;
see \cite{smirnov-werner}.
Christophe Garban kindly allowed us to include as
Appendix \ref{appgarban} his proof
for bond percolation on the square lattice,
which is partly inspired by the noise-sensitivity
arguments in~\cite{benjamini-kalai-schramm}.
It seems possible to extend the proof to a wider range of percolation models
(by changing the values of the $C_j$ random variables so they become
centered, rather than using symmetry arguments).
\begin{question}[(Estimating pivotals)]\label{qpivo}
For which models can (\ref{eqpivo}) be proved?
Can it be deduced from (\ref{eqwrsw})?
Is there a geometric argument?
\end{question}


Summing up the discussion above, we study collections of percolation
models satisfying
Assumptions \ref{assursw}.
In particular,
\textit{our results apply to
critical and near-critical
site percolation on the triangular lattice and
bond percolation on the square lattice}.

\subsection{Definition of the 
scaling limits}

Every discrete percolation configuration
(or a coloring of a tiling)
is clearly described
by the collection of quads ($=$ topological quadrilaterals) which it crosses
(and for a fixed lattice mesh in a bounded domain a finite collection
of quads would suffice).
We introduce a setup, which allows to work independently of the lattice mesh
and pass to the scaling limit,
but the idea remains the same: a percolation configuration is
encoded by a collection of quads---those which are crossed by clusters.
Our setup is inspired by the Dedekind's sections and uses the
metric and ordering on the space of quads.
The ordering comes from the obvious observation that crossing of a quad
automatically contains crossings of ``shorter'' sub-quads,
and so possible percolation configurations ($=$ quad collections)
should satisfy
certain monotonicity properties.

\begin{figure}

\includegraphics{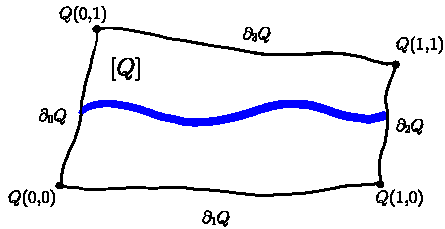}

\caption{A quad with a crossing.}
\label{figquad}
\end{figure}

We will work in a domain $D\subset\C=\R^2$.
A \textit{quad} in $D$ is a topological quadrilateral,
that is, a homeomorphism $Q\dvtx[0,1]^2 \to Q([0,1]^2)\subset D$.
We introduce some redundancy when considering
quads as parameterized quadrilaterals, but we
avoid certain technicalities.
The space of all such quads will be denoted by
$\Quad=\Quad_D$. It is a metric space under the uniform metric
\[
d(Q_1,Q_2)=\sup_{z\in[0,1]^2} |Q_1(z)-Q_2(z)|.
\]
A \textit{crossing} of $Q$ is a connected compact subset
of $[Q]:=Q([0,1]^2)$ that intersects both
opposite sides
${\partial}_0Q:=Q(\{0\}\times[0,1])$ and ${\partial
}_2Q:=Q(\{1\}\times[0,1]\})$.
We also denote the remaining two sides by
${\partial}_1Q:=Q([0,1]\times\{0\})$ and
${\partial}_3Q:=Q([0,1]\times\{1\})$,
see Figure \ref{figquad}.
%
%
The whole boundary of $Q$ we denote by ${\partial}Q$.
In the discrete setting there is no difference between
connected and path-connected crossings,
but one can imagine lattice models where
the probabilities of the two would be different in the scaling limit.
However, for the models under consideration,
every crossing in the scaling limit almost surely can be realized
by a H\"{o}lder continuous curve \cite{aizenman-burchard}, so the two notions
are essentially the same.

\begin{figure}

\includegraphics{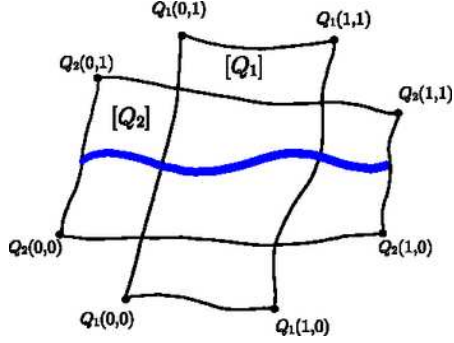}

\caption{Any crossing of $Q_2$ contains a crossing of $Q_1$, so we
write $Q_1\le Q_2$.}
\label{figorder}
\end{figure}

In addition to the structure of $\Quad_D$ as a metric space, we
will also use the following partial order on $\Quad_D$.
If $Q_1,Q_2\in\Quad_D$, we write $Q_1\le Q_2$ if
every crossing of $Q_2$ contains a crossing of $Q_1$.
The simplest example is seen in Figure \ref{figorder}.
%
%
Also, we write $Q_1<Q_2$ if there are open (in the uniform metric) neighborhoods
$U_1$ of $Q_1$ and $U_2$ of $Q_2$ in $\Quad_D$,
such that for every $Q\in U_1$ and $Q'\in U_2$
we have $Q\le Q'$.
Thus, the set of pairs $(Q,Q')\in\Quad^2$ satisfying $Q<Q'$
is the interior of $\{(Q,Q')\in\Quad^2:Q\le Q'\}$
in the product topology.

A subset $S\subset\Quad_D$ is a \textit{lower set} if,
whenever $Q\in S$ and $Q'\in\Quad_D$ satisfies $Q'< Q$, we also have
$Q'\in S$.
The collection of all closed lower subsets of~$\Quad_D$
will be denoted $\Scal_D$.
Any discrete percolation configuration $\dsamp$ in $\C$ is naturally
associated
with an element $S_\dsamp$ of $\Scal_D$: the set of all
quads for which $\dsamp$ contains a crossing
formed by an open cluster.
Thus, a percolation model on the tiling $H_\eta$ induces
a probability measure $\Pdisc_\eta$ on $\Scal_\C$
(and more generally on $\Scal_D$ for any domain $D$
contained in its domain of definition).

A topology is defined on $\Scal_D$ by specifying a subbase.
If $U\subset\Quad_D$ is topologically open and $Q\in\Quad_D$, let
\[
V_U := \{S\in\Scal_D \dvtx S\cap U\ne\varnothing\}
\]
and
\[
V^Q := \{S\in\Scal_D\dvtx Q\notin S\}.
\]
If we regard $S$ as a collection of crossed quads, then
\[
V_U := \{S \dvtx\mbox{some quads from } U \mbox{ are crossed}\}
\]
and
\[
V^Q := \{S\dvtx\mbox{quad } Q \mbox{ is not crossed}\}.
\]

Our topology of choice on $\Scal_D$ will be
the minimal topology $\mathcal T_D$ which contains every such $V_U$ and $V^Q$.
As we will see, $(\Scal_D,\mathcal T_D)$ is a metrizable compact
Hausdorff space.
The percolation scaling limit configuration will be defined as an
element of $\Scal_\C$,
and the scaling limit
measure is a Borel probability measure on~$\Scal_\C$.

Then the event of a
percolation configuration crossing a quad $Q$ corresponds
to
\[
\cross{Q} := \neg V^Q = \{S\in\Scal\dvtx Q\in S\} \subset\Scal.
\]
We will prove that if $\Quad^0\subset\Quad$ is dense in $\Quad$,
then the events $\cross{Q}$ for $Q\in\Quad^0$, generate the
Borel $\sigma$-field of $\Scal$.

\subsection{Statement of main results}

Construction of the
(subsequential) scaling limits is accomplished
for all random colorings of trivalent tilings
without additional assumptions.
Other results of this paper apply to all
\textit{percolation models on trivalent tilings},
\textit{satisfying Assumptions} \ref{assursw},
that is, the RSW estimate (\ref{eqwrsw})
and the pivotal estimate (\ref{eqpivo}).

As discussed above, we fix a collection of such percolation models
indexed by a set $\{\eta\}$,
and denote by $\Pdisc_\eta$ the (locally) discrete probability measure
on percolation configurations
for the model corresponding to parameter $\eta$.
By $|\eta|\in\R_+$, we denote the scale parameter of the model
(which enters in the RSW estimates).
Later we will show that $\Pdisc_\eta$ has subsequential scaling
limits as $|\eta|\to0$,
and will work with one of those, denoted by $\Plim$.

First, we present a simplified discrete version of the gluing theorem.
Informally, the statement is that if you fix a quad $Q_0$
and a finite length path $\alpha$ cutting it into two pieces,
then the discrete percolation configuration outside
a small neighborhood of $\alpha$ reliably predicts the
existence or nonexistence of a crossing of $Q_0$ with probability close
to one,
uniformly in the mesh size.

Following is the precise version,
estimating the random variable which is
the conditional probability of observing an open crossing
given the configuration $s$-away from $\alpha$.
As a regularity assumption, we need a constant $C(\alpha)$ such
that for every $s>0$ the set $\alpha$ can be covered by at most
$C(\alpha)/s$ balls of radius $s$.
This is the case
if and only if $\alpha$ has finite one-dimensional upper Minkowski content
$\mathcal{M}^{*1}(\alpha) = \limsup_{s\to0+}(\area\{z\dvtx\dist
(z,\alpha)<s\}) $.
%
\begin{theorem}[(Discrete gluing)]\label{tdiscrete}
Consider a collection of percolation~mo\-dels
satisfying Assumptions \ref{assursw}. 
Let $Q_0\in\Quad$ be some quad, and let $\alpha\subset[Q_0]$
be a finite union of finite length paths,
or more generally a set of finite
one-dimensional upper Minkowski content.
Let $\cross{Q_0}$ denote the event that~$Q_0$ is crossed by $\omega$, and for each $s>0$ let
${\mathcal{F}}_s$ denote the $\sigma$-field generated by the
restriction of
$\omega$ to the complement of the $s$-neighborhood of $\alpha$.
Then for every $\eps>0$,
\[
\lim_{s\searrow0} \sup_{|\eta|\in(0,s)}
\Pdisc_\eta\bigl(\eps<\Pdisc_\eta(\cross{Q_0}\mid{\mathcal
{F}}_s)<1-\eps\bigr)=0 .
\]
\end{theorem}
\begin{remark}[(Gluing is nonconstructive)]\label{rempmain}
Theorem \ref{tdiscrete} essentially states that event $\cross{Q_0}$
\textit{for any given} $\eta$
can be reconstructed with high accuracy by sampling
crossing events away from $\alpha$.
This does not lead to the noise characterization of the scaling limit,
since a priori the number of events we need to sample can grow out of control
as $|\eta|\to0$.
We address this by providing a uniform bound on the number
of events to be sampled in Proposition \ref{pmain}.
However, we do not give a recipe, and in general the gluing procedure
may still depend on $\eta$.
We deal with this possibility by assuming that there is
a~(subsequential) scaling limit, then,
since the number of possible outcomes of sampling a finite number of
crossing events is finite,
we can choose a sub-subsequence along which the same procedure is used.
The downside is that, in the current form, our results alone
cannot be used to show the uniqueness of the percolation scaling limit.
\end{remark}

The previous remark leads to the following:
\begin{question}[(Constructive gluing)]
Can one show that there is an independent of $\eta$ and preferably
constructive procedure to
predict $\cross{Q_0}$ by sampling crossing events away from $\alpha$?
\end{question}

A positive answer would provide an elegant proof of the uniqueness of
the scaling limit
(provided Cardy's formula is known), make the universality phenomenon
more accessible, and lay foundation for the renormalization theory
approach to percolation,
similar to one discussed by Langlands; see~\cite{langlands-renorm}
and the references therein.
A particularly optimistic scenario
is suggested in Conjecture \ref{conjtight}.


Theorem \ref{tdiscrete} is a toy version of our main gluing results,
and the logic behind it was folklore for a long time
(but strangely enough appeared few times in print,
as was pointed out in \cite{steif-survey}).
The proof proceeds by consequently resampling
radius $2s$ balls covering the $s$-neighborhood of $\alpha$ and using
the pivotal estimate (\ref{eqpivo}).
For this to work, the bounds have to sum up to something small,
and so with better pivotal estimates we can relax assumptions on
$\alpha$
(and the number of balls needed to cover it).
\begin{remark}[(Regularity assumptions on $\alpha$)]
When applying the pivotal estimate (\ref{eqpivo}),
we cannot relax the regularity assumption that
$\alpha$ has finite one-dimensional upper Minkowski content.
However, as discussed above, most percolation models are expected to
satisfy stronger estimates,
and then the regularity assumption can be considerably weakened
(though it cannot be dropped).
The set of pivotal points conjecturally has Hausdorff dimension~$3/4$
in the scaling limit
(proved on the triangular lattice in~\cite{smirnov-werner}),
so if the Hausdorff dimension
of the set $\alpha$ is greater than~$5/4$,
the two should intersect with positive probability
(needed correlation estimates are expected to hold).
In such a situation, there is little hope to reliably determine the crossing
without knowledge of $\alpha$. 
To be more rigorous, instead of pivotals one has to work with the
percolation spectrum,
but it has the same dimension (see~\cite{garban-pete-schramm} for the
definition and discussion),
so $\alpha$ has to be of dimension at most~$5/4$.
On the other hand, our proof works for sets $\alpha$ of dimension
smaller than~$2$ minus dimension of the pivotal points,
provided also that $\alpha\cap{\partial}Q_0$ is finite.
So~$5/4$ seems to be the correct-dimensional threshold,
and this can be established for
the critical site percolation on the triangular lattice.
\end{remark}
\begin{remark}[(Crossings with finite intersection
number)]\label{rem23p23}
Suppose that a smooth curve $\alpha$ cuts the quad $Q_0$ into two
quads $Q_+$ and $Q_-$.
The following informal discussion concerns the
critical percolation scaling limit,
assuming its existence and conformal invariance
(but it can be easily made more precise by working
with the limits of discrete configurations and events).
First observe that open clusters inside $Q_\pm$
intersect the smooth boundary $\alpha$ on a set of Hausdorff dimension $2/3$
(and we have a control for the correlation of points inside it).
Being independent, two large open clusters on opposite sides
of $\alpha$ have a positive probability of touching.
But if they touch at some point $x\in\alpha$,
we cannot have a large closed cluster passing through $x$ and
separating them.
Indeed, that would mean a four-arm event occurring at the
\textit{pivotal} point $x$ on the smooth curve $\alpha$, and since $a_4=5/4$,
the set of pivotal points has Hausdorff dimension $3/4$ and so
almost surely does not intersect the smooth $\alpha$.
Thus, the two open clusters in $Q_\pm$ which touch on~$\alpha$
are in fact parts of the same large open cluster in $Q_0$.
We conclude, that given a~smooth curve $\alpha$,
with positive probability there is a crossing of $Q_0$ which intersects
$\alpha$
only \textit{once}.
\end{remark}

The last remark says that for a smooth curve $\alpha$,
if we condition on the quad being crossed, the probability
of finding a crossing intersecting $\alpha$ only finitely many times
is positive.
We can ask whether this conditional probability is one, leading to the
following:
\begin{conjecture}[(Finite intersection property)]\label{conjfinite}
Let $Q$ be a quad and~$\alpha$ a smooth curve.
Then in the scaling limit, given existence of a crossing,
there is almost surely a crossing intersecting $\alpha$ only a finite
number of times.
\end{conjecture}

Alternatively, one can formulate a discrete version, which is a priori stronger
(for an equivalent one we would have to count small scale
packets of intersections):
\begin{conjecture}[(Tightness of intersections)]\label{conjtight}
Let $Q$ be a quad and~$\alpha$ a~smooth curve.
For a discrete percolation measure $\Pdisc_\eta$, let the random variable
$N_\eta=N_\eta(Q,\alpha)$
be the minimal possible number of intersections with~$\alpha$ for
a~crossing of $Q$ (and $0$ if $Q$ is not crossed).
Then $N_\eta$ is tight as $|\eta|\to0$.
\end{conjecture}

Proving either of these conjectures would not only simplify our proof
of the main theorem,
but also provide an easy and constructive ``gluing procedure'':
we take all ``excursions'' ($=$ crossings in the complements of~$\alpha$) and check whether
a crossing of $Q$ can be assembled from finitely many of them
(as discussed in Remark \ref{rem23p23},
the finite number of excursions will necessarily connect into one crossing).
Such a procedure would give a very short proof of the uniqueness of the
full percolation scaling limit,
provide an elegant setup for the renormalization theory of percolation
crossings,
and show that the ``finite model'' discussed by Langlands and Lafortune
in~\cite{langlands-finite}
indeed approximates the critical percolation.
Note that if conjectures above fail, the constructive procedure for
gluing crossings will have to deal with countably many excursions
intersecting $\alpha$ on a Cantor set $C$.
Then the criterion for possibility of gluing the excursions into one crossing
is likely to be in terms of some capacitary (or dimensional)
characteristics of $C$.

The critical percolation scaling limit is conjectured to be conformally
invariant,
which was established for site percolation on the triangular lattice~\cite{smirnov-perc,smirnov-cras}.
In Remark \ref{rem23p23},
we mentioned that, for a given quad $Q_0$,
a~smooth curve $\alpha$ has
the property that two large clusters inside two components of
$[Q_0]\setminus\alpha$
touch each other (on $\alpha$)
with positive probability.
One can ask, whether an arbitrary curve $\alpha$ would always enjoy
this property,
and the answer is negative.
Moreover, under conformal invariance assumption,
this property (equivalent to positive probability of finding a crossing
intersecting $\alpha$ only finitely many times)
can be reformulated in terms of multifractal properties of
harmonic measure on two sides of~$\alpha$.
Below we present an informal argument to this effect,
note also similar discussions in the SLE context in \cite
{garban-rohde-schramm}.

To simplify the matters, we discuss a related (and probably equivalent) property
that the expected number of such touching points
for two large clusters inside two components of $[Q_0]\setminus\alpha$
is positive (or even infinite).
Let $f(a_+,a_-)$ be the two-sided dimension spectrum of harmonic
measure on $\alpha$; see \cite{astala-prause-smirnov}.
Roughly speaking, $f(a_+,a_-)$ is the dimension of the set of $z\in
\alpha$ where harmonic measures $\omega_\pm$
on two sides of $\alpha$ have power laws $a_\pm$:
\[
\omega_+(B(z,r)) \approx r^{a_+} ,\qquad \omega_-(B(z,r)) \approx r^{a_-} .
\]
Consider a set $E\subset\alpha$ such that harmonic measures $\omega
_\pm$
have power laws~$a_\pm$ on~$E$.
Cover most of $E$ by a collection of balls $B_j$ of small radii $r_j$.
By conformal invariance,
the probability of a ball $B_j$ to be touched
by a large open cluster on the ``$+$ side'' of $\alpha$ is
governed by its harmonic measure and Cardy's formula \cite{smirnov-cras}
that gives the exponent $1/3$:
\[
P_+(B_j) \approx\omega_+(B_j)^{1/3} \approx r_j^{a_+/3} .
\]
Combining with the similar estimate on the ``$-$ side,'' we write
an estimate for the expected number of balls $B_j$
touched by large clusters on both sides:
\[
\sum_j P_+(B_j)P_-(B_j) \approx\sum_j r_j^{(a_++a_-)/3} .
\]
Here the right-hand side is small (large) when
$(a_++a_-)/3$ is bigger (smaller) than $\dim(E)$,
since $\sum_j r_j^{\dim(E)}\approx1$.
We conclude that the expected number of points
touched by large clusters on both sides
is positive if and only if we can find $E$
such that $(a_++a_-)/3<\dim(E)$.
\begin{remark}[(Sharpness of the finite intersection property)]
It seems that positive probability of having a crossing intersecting
$\alpha$ a finite number of times
is roughly equivalent to the existence of positive exponents $a_+$ and~$a_-$
such that the two-sided multifractal spectrum of harmonic measure on
$\alpha$
satisfies $3 f(a_+,a_-)>a_++a_-$; see \cite{astala-prause-smirnov}
for a discussion of such spectra.
The latter property seems to fail even for some nonsmooth $\alpha$
of dimension close to one,
so we do not expect the factorization property
to be equivalent to the finite intersection property,
but rather to be a weaker one.\vspace*{-2pt}
\end{remark}

Before proving the main theorem, we discuss a new framework for the
percolation scaling limit
and show that collection of discrete random colorings is precompact,
so in this framework the subsequential scaling limits exist for any
random model.

Recall that we work with $\Scal_D$,
which is the collection of all closed lower
(i.e., monotonicity of the crossing events is obeyed)
subsets of $\Quad_D$,
which one should think of as sets
of quads crossed.
The topology $\mathcal T_D$ is the minimal one containing
all sets $V_U$ (of $S$ such that some quads in $U$ are crossed)
and~$V^Q$ (of $S$ such that $Q$ is not crossed).

We now list some important properties of $\Scal_D$.\vspace*{-2pt}
\begin{theorem}[(The space of percolation configurations)]\label{tpscal}
Let $D\subset\hat\C$ be topologically open and nonempty.
\begin{longlist}[(2)]
\item[(1)] $(\Scal_D,\mathcal T_D)$ is a compact metrizable Hausdorff space.
%
\item[(2)] Let $A$ be a dense subset of $\Quad_D$.
If $S_1,S_2\in\Scal_D$ satisfy $S_1\cap A=S_2\cap A$, then $S_1=S_2$.
Moreover, the $\sigma$-field generated by $V^Q$, $Q\in A$ is
the Borel $\sigma$-field of $(\Scal_D,\mathcal T_D)$.
\item[(3)] If $S_1,S_2,\ldots$ is a sequence in $\Scal_D$ and $S\in
\Scal_D$, then
$S_j\to S$ in $\mathcal T_D$ is equivalent to
$S=\limsup_j S_j=\liminf_j S_j$.
\item[(4)] Let $D'\subset D$ be a subdomain, and for $S\in\Scal_D$
let $S':= S\cap\Quad_{D'}$. Then $S\mapsto S'$ is a continuous map from
$\Scal_D$ to $\Scal_{D'}$.
\end{longlist}
\vspace*{-2pt}
\end{theorem}

In the above, $\liminf_j S_j$ is the set of $Q\in\Quad_D$ such that
every neighborhood of $Q$ intersects all but finitely many of the
sets $S_j$, and $\limsup_j S_j$ is the set of $Q\in\Quad_D$
such that every neighborhood of $Q$ intersects infinitely many of the
sets~$S_j$.\vspace*{-2pt}
\begin{remark}[(Applicability of Theorem \ref{tpscal})]\label
{remconditionsforscaling}
The result concerns the topological properties of the space of
percolation configurations (colorings of tilings)
and involves no probability measures.
Thus, it can be applied to \textit{any random coloring} of a trivalent tiling
(into open and closed tiles),
yielding results for dependent models, like
the Fortuin--Kasteleyn percolation.\vspace*{-2pt}
\end{remark}

It follows from the theorem above that the space of continuous
functions on $(\Scal_D,\mathcal T_D)$ is separable,
and hence the unit ball in its dual space
$M(\Scal_D)$ of Borel measures is compact\vadjust{\eject}
(by the Banach--Alaoglu theorem, Theorem~3.15 in \cite{rudin-fan}),
metrizable
(by Theorem 3.16 in \cite{rudin-fan})
and obviously Hausdorff in the weak-$*$ topology.
The space of probability measures being a closed convex subset of the
unit ball, we arrive at the following.
\begin{corollary}[(The space of percolation measures)]\label{cormeasures}
The space\break $\operatorname{Prob}(\Scal_D)$
of Borel probability measures on $(\Scal_D,\mathcal T_D)$
with weak-$*$ topology
is a compact metrizable Hausdorff space.
\end{corollary}

Discrete percolation measures $\Pdisc_\eta$ clearly
belong to $\operatorname{Prob}(\Scal_D)$
and hence the existence of percolation
subsequential scaling limits
(or, more generally, for any random colorings)
is an immediate corollary:
\begin{corollary}[(Precompactness)]\label{corsubseqexist}
The collection of all discrete measures $\{\Pdisc_\eta\}$
(corresponding to random colorings of trivalent tilings)
is precompact in the topology of weak-$*$ convergence
on $\operatorname{Prob}(\Scal_D)$.
\end{corollary}
\begin{remark}[(Nets vs. sequences)]\label{remseqnet}
We can make $\{\eta\}$ a directed set
by writing $\eta_1\preceq\eta_2$
whenever $|\eta_1|\ge|\eta_2|$.
Then the scaling limit $\lim_{|\eta|\to0}\Pdisc_\eta$ is
just the limit of the net $\Pdisc_\eta$.
Since the target space $\operatorname{Prob}(\Scal_D)$
is metrizable and hence first-countable,
it is enough to work with the notion of
convergence along sequences, rather than nets.
In particular, if $\lim_{j\to\infty}\Pdisc_{\eta_j}=\Plim$
for every sequence $\eta_j$ with $\lim_{j\to\infty}|\eta_j|=0$,
then $\lim_{|\eta|\to0}\Pdisc_{\eta}=\Plim$.
\end{remark}
\begin{remark}
Our construction keeps information about all macroscopic percolation events,
and so gives the \textit{full percolation scaling limit}
at or near criticality.
For the critical site percolation on triangular lattice,
Garban, Pete and Schramm explain in Section 2 of \cite
{garban-pete-schramm-pivotal}
that the scaling limit is unique, appealing to the work \cite
{camia-newman-cmp} of Camia and Newman.
Basically, we know the probabilities of individual crossing events, and
it allows to reconstruct the whole picture, albeit in a difficult way.
It would be interesting to establish the uniqueness using an
appropriate modification of Proposition~\ref{pmain}.\looseness=1
\end{remark}

So that our theorems apply even to models where the uniqueness of the
scaling limit has not been proved yet,
we work with one of the subsequential limits~$\Plim$---a~Borel
probability measure on $\Scal_\C$, provided by Corollary~\ref{corsubseqexist}.
By Theorem~\ref{tpscal}, the Borel $\sigma$-field of $\Scal_D$ is
$\sigma(\cross{Q}\dvtx Q\in\Quad_D)$,
that is, generated by the crossing events inside $D$.
Let ${\mathcal{F}}_D:=\sigma(\cross Q\dvtx Q\in\Quad_D)$ also denote the
corresponding
subfield of the Borel $\sigma$-field of $\Scal_\C$.
\begin{theorem}[(Factorization)]\label{tfactor}
Consider a collection of percolation mo\-dels,
satisfying Assumptions \ref{assursw}. 
Let $D$ be a domain, and let $\alpha\subset\C$ be a~finite union of
finite length paths with finitely many double points.\vadjust{\eject}
Denote the components of $D\setminus\alpha$ by $D_j$.
Then for any subsequential scaling limit
%
%
\begin{equation}
\label{eqprodcont}
{\mathcal{F}}_D = {\mathcal{F}}_{D\setminus\alpha} = \bigvee
_j{\mathcal{F}}_{D_j} ,
\end{equation}
up to sets of measure zero.\vspace*{-2pt}
\end{theorem}

Again, we note that this does not hold when $\alpha$ is a sufficiently
wild path.

Fix a subsequential scaling limit of the critical percolation for
the bond model on the square lattice or
for the site model on the triangular lattice.
By the results above for a domain $D$,
it is described by a probability space
$(\Scal_D, {\mathcal{F}}_D, \Pdisc_D)$,
which is invariant under translations of the plane
(in fact, under all conformal transformations, but it is so far established
for the triangular lattice only),
depends continuously on $D$ and satisfies (\ref{eqprodcont}).
Then in the language of Tsirelson
(see \cite{tsirelson-nonclassical}, Definition 3d1)
it is a \textit{noise}
or a \textit{homogeneous continuous product of probability spaces}
with a Boolean base given by an appropriate algebra
of piecewise-smooth planar domains
(e.g., generated by rectangles).

There is a well developed and beautiful theory of noises,
and we refer the reader to two expositions
\cite{tsirelson-nonclassical,tsirelson-lnm} by Tsirelson.
In particular, every noise~can be decomposed
into a ``classical,''~or ``stable'' part
and a ``nonclassical,'' or ``sensitive'' part.
White noise is purely classical,
and purely nonclassical noises are called ``black'' by Tsirelson.
The black noises are difficult to construct,
with the first example given by Tsirelson and Vershik in \cite
{tsirelson-vershik}.
As pointed out by Tsirelson in
\cite{tsirelson-experiment} and in Remark 8a2 of \cite{tsirelson-lnm},
percolation noise has to be nonclassical, and so our paper
provides the first genuinely two-dimensional example of a black
noise.\vspace*{-2pt}
\begin{corollary}[(Percolation is a noise)]
Thus, we conclude that,
in the language of Tsirelson \cite{tsirelson-nonclassical},
any subsequential scaling limit as above is a \textit{noise}
with a Boolean base given by an appropriate algebra
of piecewise-smooth planar domains
(e.g., generated by rectangles).
Therefore, it has to be a
\textit{black noise}, as explained in \cite{tsirelson-lnm}, Remark
\textup{8a2}.\vspace*{-2pt}
\end{corollary}

By the conjectured universality, all
the critical percolation models should have the same scaling limit
and so correspond to the same black noise.
The situation with near-critical noises is less clear,
so the following question was proposed by Christophe Garban
(the answer is expected to be negative):\vspace*{-2pt}
\begin{question}\label{quniv}
Are the noises arising from the critical and near-critical models isomorphic,
or do we get a family of different noises?\vspace*{-2pt}
\end{question}

The notion of isomorphism for noises is discussed by Tsirelson
in Section~4a of~\cite{tsirelson-nonclassical}.
Roughly speaking, we ask
\textit{whether one can get a near-critical noise from the critical one
by a local deterministic procedure}?
Note that the current construction of the near-critical percolation
scaling limit
would use the critical percolation along with some extra randomness,
as in \cite{camia-fontes-newman-loops,garban-pete-schramm-icmp}.\vadjust{\eject}

Our construction of the scaling limit applies to general random colorings,
and much of what we do afterwards mostly uses ``arms estimates,''
which are known for a wider class of models, than just percolation.
Therefore, it is logical to ask:
\begin{question}\label{qdependent}
$\!\!$To what extent out results apply to the models
with~de\-pendence?
One has to be more careful in formulating the results,
since now the configuration inside a domain depends on the boundary conditions,
but this can be addressed by working with a full-plane model and stripe domains.
\end{question}

\section{The discrete gluing theorem}\label{secdiscretegluing}

In this section, we prove Theorem \ref{tdiscrete}.

Applying several times in succession Lemma \ref{lcont},
we can find a smooth (i.e., given by a diffeomorphism) quad $Q$,
such that the crossing events for $Q_0$ and $Q$ are arbitrarily close.
By the same lemma, as $\varepsilon\to0$, the crossing event for $Q$
is well approximated by the crossing event of its
perturbation $Q^\varepsilon:=Q([\varepsilon,1-\varepsilon]^2)$.
On the other hand, $\alpha$ has finite length
and so for almost every~$\varepsilon$
the intersection $\alpha\cap{\partial}Q^\varepsilon$ is finite.
This can be deduced, for example, from the area theorem in the
geometric measure theory
\cite{federer-book,mattila-book}, which implies that an orthogonal
projection of a
rectifiable curve on a given line covers almost every point at most
finitely many times.

Thus, approximating if necessary, we can assume that
$\alpha$ intersects ${\partial}Q_0$ at finitely many points
$x_i$, $i=1,\ldots,k$.
Fix $\eps_0>0$, then
by the RSW estimate~(\ref{eqwrsw}), there is an $r>0$
such that, if $0<|\eta|<r$, then
the probability under~$\Pdisc_\eta$
that there is a crossing of $Q_0$
which intersects any of the disks~$B(x_i,r)$
is less than $\eps_0$. Fix such an $r$,
remove from the curve $\alpha$ the disks,
the rest will be away from the boundary,
simplifying future estimates:
\[
\alpha' := \alpha\bigm\backslash\bigcup_i B(x_i,r)
\]
and
\[
d := \inf\{|x-y|\dvtx x\in\alpha', y\in{\partial}Q_0\} > 0 .
\]

\begin{figure}

\includegraphics{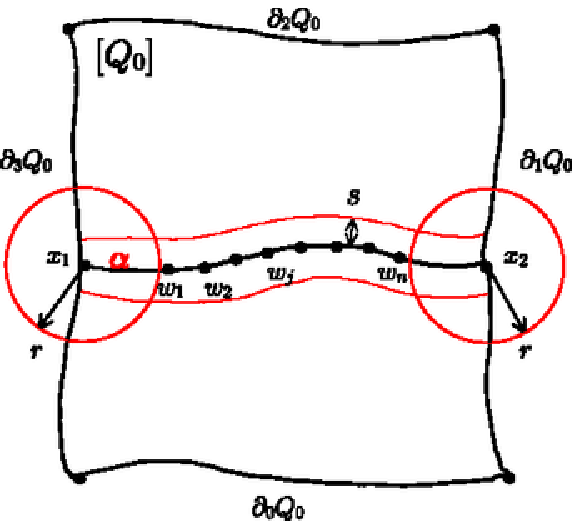}

\vspace*{-2pt}
\caption{Setup for the proof of the discrete gluing theorem in the
simplest case.}
\label{figdiscrete-setup}
\vspace*{-2pt}
\end{figure}

Fix $s\in(0,d/4)$.
By the regularity assumption, we can choose $n\le C(\alpha)/s$ balls
of radius $s$,
entirely covering $\alpha$. Denote their centers by $w_1,w_2,\ldots,w_n$,
and set $B_j:=B(w_j,2s)$; see Figure \ref{figdiscrete-setup}. 

%

Let $M$ denote the
union of the set of tiles of $H_\eta$
that intersect the $s$-neighborhood of $\alpha'$ minus $\bigcup_i B(x_i,r)$.
Assume $|\eta|<s$, so that $M$ would be contained in the
$2s$-neighborhood of $\alpha'$.
Let~$\omega$ and $\omega'$ be samples from~$\Pdisc_\eta$ that agree
on the complement
of $M$ and are independent in $M$.
For $j=0,1,\ldots,n$, let $\omega_j$ be the configuration that agrees
with $\omega'$ on the tiles of~$H_\eta$ lying inside $B_1\cup B_2\cup
\cdots\cup B_j$
and agrees with~$\omega$ elsewhere.
Then each~$\omega_j$ is a fair sample from~$\Pdisc_\eta$.
Also observe that $\omega_0=\omega$ and $\omega_n$ differs from
$\omega'$ only
inside $\bigcup_i B(x_i,r)$.

\begin{figure}[b]
\vspace*{-2pt}
\includegraphics{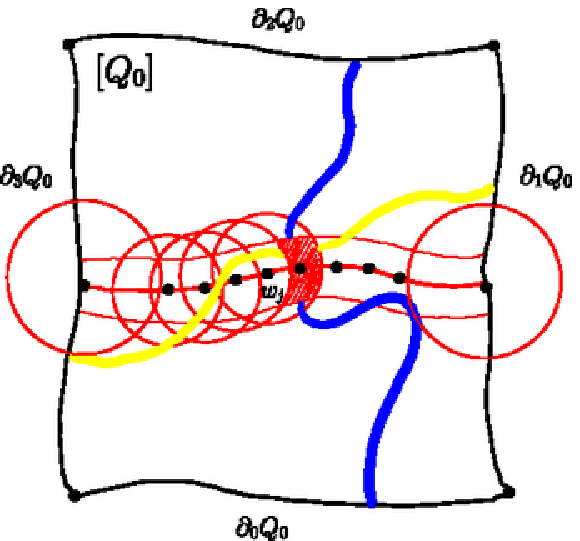}
\vspace*{-2pt}
\caption{If resampling a part of the ball $B_j$ changes the crossing
event, than there are
two open crossings from $B_j$ to ${\partial}_0Q_0$ and ${\partial}_2Q_0$,
and two dual closed crossings from $B_j$ to ${\partial}_0Q_0$ and
${\partial}_2Q_0$, a~four arm event.}
\label{figdiscrete-pivotal}
\end{figure}

Our goal is to estimate
\[
\mathbf{P}[\omega\in\cross{Q_0},\omega'\notin\cross{Q_0}],
\]
which we will do by successively comparing $\omega_{j-1}$ to $\omega_{j}$.
Fix $j\in\{1,2,\ldots,n\}$.
In order for $\{\omega_{j-1}\in\cross{Q_0}\}\cap\{\omega_{j}\notin
\cross{Q_0}\}$ to hold,
there must be a crossing in $\omega_{j-1}$ from ${\partial}_0Q_0$ to
${\partial}_2 Q_0$ that goes through
$B_j$ and there must be a dual closed crossing in $\omega_j$ from
${\partial}_1Q_0$ to
${\partial}_3 Q_0$ that goes through $B_j$; see Figure \ref
{figdiscrete-pivotal}.
%
%
Thus, in $\omega_j$ we have the
four arm event from ${\partial}B_j$ to ${\partial}{Q_0}$.
Since the radius of $B_j$ is $2 s$ and the distance from $B_j$ to
${\partial}{Q_0}$ is
at least $d-2 s>d/2$,
by assumption (\ref{eqpivo}) we have
\begin{eqnarray*}
\mathbf{P}[\omega_{j-1}\in\cross{Q_0},\omega_{j}\notin\cross
{Q_0}]
&\le&
\prsw^{4}_{\eta}(z,2s,d/2)
< \biggl(\frac{2s}{d/2}\biggr) \cdot\rsww_4(2s,d/2).
\end{eqnarray*}
Summing over $j$, we conclude
that
\begin{eqnarray*}
\mathbf{P}[\omega_0\in\cross{Q_0}, \omega_n\notin\cross{Q_0}]
&\le& \sum_{j=1}^n
\mathbf{P}[\omega_{j-1}\in\cross{Q_0},\omega_{j}\notin\cross
{Q_0}]\\
&\le& n \cdot\biggl(\frac{2s}{d/2}\biggr) \cdot\rsww_4(2s,d/2)\\
&\le& {\frac{C(\alpha)}s} \cdot{\frac{4s}{d}} \cdot\rsww
_4(2s,d/2)\\
&\le& {\frac{4C(\alpha)}{d}} \cdot\rsww_4(2s,d/2) .
\end{eqnarray*}
%
By Assumption \ref{assursw} the right-hand side
tends to $0$ with $s$, and so we can make sure that
$ \mathbf{P}[\omega_{0}\in\cross{Q_0},\omega_{n}\notin\cross
{Q_0}]<\eps_0$
by choosing $s$ sufficiently small.
This is the only place in the paper, where we use
the pivotal estimate (\ref{eqpivo}).

Now, recall that $\omega_0=\omega$, and $\omega_n$ differs from
$\omega'$ only near endpoints of curves in $\alpha$,
that is, inside $\bigcup_i B(x_i,r)$.
By our choice of $r$, we have the estimate
$\mathbf{P}[\omega_n\in\cross{Q_0},\omega'\notin\cross
{Q_0}]<\eps_0$.
Summing the above, we get the bound
\[
\mathbf{P}[\omega\in\cross{Q_0},\omega'\notin\cross{Q_0}]<2 \eps_0.
\]
Now note that by first conditioning
on ${\mathcal{F}}_s$ and then using the fact that $\omega$ and
$\omega'$ are conditionally independent given ${\mathcal{F}}_s$, we obtain
\[
2 \eps_0>\mathbf{P}[\omega\in\cross{Q_0},\omega'\notin\cross{Q_0}]=
\E\bigl[\Pdisc_\eta(\cross{Q_0}\mid{\mathcal{F}}_s) \bigl(1- \Pdisc_\eta
(\cross{Q_0}\mid{\mathcal{F}}_s)\bigr)\bigr].
\]
In particular, by Chebyshev inequality,
\[
\Pdisc_\eta\bigl(\eps<\Pdisc_\eta(\cross{Q_0}\mid{\mathcal
{F}}_s)<1-\eps\bigr)
<2 \eps_0/\bigl(\eps(1-\eps)\bigr).
\]
Since $\eps_0$ was an arbitrary positive number, this completes the proof.

\vspace*{-2pt}\section{The space of lower sets}
\label{sechered}

The main goal of this section is to provide an abstract setup for the
percolation
scaling limit and to prove Theorem~\ref{tpscal}.
We discuss topology on the space $\Scal_X$ of closed lower
subsets of an abstract partially ordered topological space $X$.
In the percolation case, $X$ corresponds to the space $\Quad_D$ of quads,
and the space $(\Scal_X,\topo_\Scal)$ to the space
$(\Scal_D,\topo_D)$ of percolation configurations
in a domain $D\subset\C$.

With these conventions,
Theorem \ref{tpscal} directly follows from Theorem \ref{tscal} below.
The latter applies to the space of percolation configurations
(colorings of tilings),
since $\Quad_D$ is clearly second-countable,
(\ref{eopen}) holds by definition and~(\ref{esb}) holds
since a quad $Q_0$ can be easily approximated by smaller quads
[e.g., by quads $Q^q$ with $q=(-\epsilon,\epsilon,1+\epsilon
,1-\epsilon)$ from the proof of Lemma \ref{lcrossingOpen}].

We now start the (mostly classical)
abstract construction, inspired by the Dedekind's sections.
Let $(X,\tau)$ be a topological space, let $\clo_X$ denote
the collection of topologically closed subsets of $X$.\vadjust{\eject}
For a topologically open $U\subset X$ and a compact $K\subset X$, let
\begin{eqnarray*}
W_U &:=& \{F\in\clo_X\dvtx F\cap U\ne\varnothing\},\\
W^K &:=& \{F\in\clo_X\dvtx F\cap K=\varnothing\}.
\end{eqnarray*}
Let $\hat\topo$ be the topology on $\clo_X$ generated by all such
sets $W^K$ and $W_U$.
It is a more convenient version of the Vietoris topology,
called \textit{Fell topology}
or \textit{topology of closed convergence}.

The following result is Lemma 1 in \cite{fell-topology}.
\begin{proposition}\label{pcompact}
$(\clo_X,\hat\topo)$ is compact.
\end{proposition}

For completness, we include a short proof, based
on Alexander's subbase theorem (see, e.g., Theorem 5.6 in \cite{kelley-book}),
which states that if $\mathcal B$ is a subbase for the topology
on a space $X$ and every cover of $X$ by elements of $\mathcal B$
contains a finite subcover, then $X$ is compact.
\begin{pf*}{Proof of Proposition \ref{pcompact}}
By Alexander's subbase theorem, it suffices to show that every cover
of $\clo_X$ by sets of the form $W_U$ (where $U\in\tau$)
and $W^K$ (where $K\subset X$ is compact) has a finite subcover.
Let
\[
\{W_U\dvtx U\in I\}\cup\{W^K\dvtx K\in J\} \supset\clo_X
\]
be such a cover.
Observe that $\bigcup_{U\in I} W_U= W_{\tilde U}$, where $\tilde
U=\bigcup_{U\in I} U$.
Let $F:=X\setminus\tilde U$, and note that $F$ is topologically closed.
Since $F\in\clo_X\subset W_{\tilde U}\cup\bigcup_{K\in J} W^K$,
we know that $F\in\bigcup_{K\in J} W^K$, namely, there is some
$K_0\in J$ with
$F\in W^{K_0}$. In other words, $F\cap{K_0}=\varnothing$; that is,
${K_0}\subset\tilde U$.
Since ${K_0}$ is compact and $\{U\dvtx U\in I\}$ covers ${K_0}$,
there is some finite $I'\subset I$ such that $\{U\dvtx U\in I'\}$
covers ${K_0}$.
Then $\{W_U\dvtx U\in I'\}\cup\{W^{K_0}\}$ is a finite cover of $\clo_X$.
\end{pf*}

Now suppose that ``$<$'' denotes a partial order on $X$ such that
%
%
\begin{equation}\label{eopen}
\mbox{
$\mathcal R:=\{(x,y)\in X^2\dvtx x<y\}$ is a topologically open subset
of $X^2$.}
\end{equation}
Define $\mathcal R^x:=\{y\,{\in}\,X\dvtx y<x\}$ and $\mathcal
R_x:=\{y\,{\in}\,X\dvtx x<y\}$.
Then for every \mbox{$x\,{\in}\,X$} these sets are topologically open.
A set $H\subset X$ is said to be \textit{lower}
if $x\in H$ implies $\mathcal R^x\subset H$.
Let $\Scal_X$ denote the collection of all topologically closed lower sets.
\begin{proposition}\label{plower}
Assuming\vspace*{2pt} (\ref{eopen}),
$\Scal_X\subset\clo_X$ is closed with respect to the topology $\hat
\topo$;
that is, $\clo_X\setminus\Scal_X \in\hat\topo$.
\end{proposition}
\begin{pf}
Set $U_x:= W^{\{x\}}\cap W_{\mathcal R_x}$ and
\[
U := \bigcup_{x\in X}U_x=\{F\dvtx\exists x,y \mbox{ with } x<y, y\in F,
x\notin F\}.
\]
Then $U$ is topologically open and $\Scal_X=\clo_X\setminus U$.\vadjust{\eject}
\end{pf}

Let $\topo_\Scal$ be the topology of
$\Scal_X$ as a subspace of $\clo_X$.
It is the topology generated by the sets
$V_U:=W_U\cap\Scal_X$ and $V^K:=W^K\cap\Scal_X$,
for topologically open $U$ and compact $K$ subsets of $X$.
As we will later see, in our setting it is sufficient to restrict ourselves
to one-point compact sets, working with~$V^{\{x\}}$,
which we will often abbreviate $V^x$.

For simplicity, we will also abbreviate
$\Scal:=\Scal_X$ and $\topo:=\topo_\Scal$.
Propositions~\ref{pcompact} and~\ref{plower} imply that
$(\Scal,\topo)$ is compact.

In general, $\Scal$ does not have to be
a Hausdorff space, even if $X$ is Hausdorff.
Indeed, the sufficient condition would be
for $X$ to be Hausdorff and compact,
but for the application we have in mind, $X$ is not compact.
However, the following result gives another sufficient condition
for $\Scal$ to be Hausdorff.\vspace*{2pt}
\begin{proposition}\label{pHau}
Suppose, in addition to (\ref{eopen}), that
%
%
\begin{equation}\label{esb}
\forall{x\in X}\qquad
x\in\closure{\mathcal R}{}^x .
\end{equation}
Then $(\Scal,\topo)$ is a Hausdorff space.
\end{proposition}
\begin{pf}
Let $H_1,H_2\in\Scal$ be different, then
one of the two differences $H_1\setminus H_2$ and $H_2\setminus H_1$ is
nonempty,
without loss of generality the first one.
Take some $x\in H_1\setminus H_2$.
Since $X\setminus H_2$ is topologically open in $X$ and contains $x\in
\closure{\mathcal R}{}^x$,
we can find some
$y\in\mathcal R^x\cap(X\setminus H_2)=\mathcal R^x\setminus H_2$.
Then $V_{\mathcal R_y}$ is disjoint from $V^{y}$.
Moreover, $H_1\in V_{\mathcal R_y}$ and $H_2\in V^{y}$.
So there are disjoint topologically open sets
in $\topo$ containing $H_1$ and $H_2$, respectively,
and the space is Hausdorff.\vspace*{2pt}
\end{pf}
\begin{remark}\label{raltsubbase}
Let $\topo'$ be the topology on $\clo_X$ generated by
the sets $W_U$, $U\in\topo_X$ and
$W^{\{x\}}$, $x\in X$. Since this topology is
coarser than $\hat\topo$, Proposition~\ref{pcompact}
implies that $(\clo_X,\topo')$ is compact.
The proof of Proposition \ref{plower} shows that
$\Scal$ is a topologically closed subset of $\clo_X$, also with
respect to~$\topo'$. Moreover, assuming~(\ref{esb}),
the proof of Proposition \ref{pHau}
shows that $\Scal$ is Hausdorff with respect to~$\topo'$.
It is well known
(and not hard to verify; see Theorem~5.8 in~\cite{kelley-book})
that if a compact topology on a space is finer than a~Hausdorff
topology on the same space, then they must be equal.
We conclude that the topology on $\Scal$
induced by $\hat\topo$ is the same as that induced
by~$\topo'$ whenever (\ref{esb}) and~(\ref{eopen}) hold.\vspace*{2pt}
\end{remark}
\begin{lemma}\label{ldensedetermines}
Suppose that (\ref{eopen}) and (\ref{esb}) hold,
and that $X_0$ is a dense subset of $X$.
If $H_1,H_2\in\Scal$ and $H_1\ne H_2$, then
$H_1\cap X_0\ne H_2\cap X_0$.
\end{lemma}
\begin{pf}
Take $x\in H_1\setminus H_2$, if the latter is nonempty.
Then (\ref{esb}) implies that $\mathcal R^x\setminus H_2\ne
\varnothing$.
Since $X_0$ is dense and $\mathcal R^x\setminus H_2$
is topologically open and nonempty, we have
$X_0\cap\mathcal R^x\setminus H_2\ne\varnothing$.
Because $\mathcal R^x\subset H_1$, this implies that
$H_1\cap X_0\ne H_2\cap X_0$. A symmetric
argument applies if $H_2\setminus H_1\ne\varnothing$.
\end{pf}

In the next lemma, we characterize convergence of nets in $\Scal_X$.
As we shall see, $(\Scal_X,\topo_\Scal)$ is first countable,
thus convergence of sequences actually leads to the same
properties, but we decided to work with the more general notion for now,
since it does not lead to additional difficulties.

Suppose that ${\mathcal D}$ is a directed set
(i.e., partially ordered by ``$\preceq$''
so that every two elements have a common upper bound),
and $\{Y_n\}_{n\in{\mathcal D}}$ is a net (i.e., a sequence indexed
by ${\mathcal D}$) of subsets of $X$.
We write $\limsup_n Y_n$ for the set of all
$x\in X$ with the property that for every topologically open
$U\subset X$ containing $x$ and for every $m\in{\mathcal D}$ there
is some $n\succeq m$ in ${\mathcal D}$ satisfying $Y_n\cap U\ne
\varnothing$.
Similarly, $\liminf_n Y_n$ is the set of all $x\in X$ with
the property that for every topologically open $U\subset X$ containing $x$
there is some \mbox{$m\in{\mathcal D}$} such that $Y_n\cap U\ne\varnothing$
for every $n\succeq m$. We write $Y=\lim_n Y_n$ if
$Y=\limsup_n Y_n=\liminf_n Y_n$.\vspace*{2pt}
\begin{lemma}\label{llim}
Assume that (\ref{eopen}) and (\ref{esb}) hold,
and that $S_n$, $n\in{\mathcal D}$, is a net
with $S_n\in\Scal$ for every $n$.
\begin{longlist}[(2)]
\item[(1)]
If $\lim_n S_n\subset X$ exists, then $\lim_n S_n\in\Scal$.
\item[(2)]
If the net $S_n$ converges to $S$ with respect to $\topo$,
then $S=\lim_n S_n$.
\item[(3)] Conversely, if $S=\lim_n S_n$,
then the net $S_n$ converges to $S$ with respect to~$\topo$.
\end{longlist}
\end{lemma}
\begin{pf}
To prove (1), suppose that $S=\lim_n S_n$, $x\in S$ and $y<x$.
Then \mbox{$x\in\mathcal R_y$}. Therefore, there is some $m\in{\mathcal D}$ such
that $S_n\cap\mathcal R_y\ne\varnothing$ for every $n\succeq m$.
Since $S_n\in\Scal$, this implies that $y\in S_n$ for every $n\succeq m$.
Thus, $y\in S$, which proves that $S$ is a lower set.
The definition of $\lim_n S_n$ makes it clear that $S$ is
topologically closed.
This proves (1).

To prove (2), suppose that $S_n$ converges to $S$ with respect to
$\topo$.
Let $x\in S$. If $U\subset X$ is topologically open and contains $x$, then
$S\in V_U$. Therefore, there is some $m\in{\mathcal D}$ such that
$S_n\in V_U$ for every $n\succeq m$. But $S_n\in V_U$ is equivalent
to $S_n\cap U\ne\varnothing$. Therefore, $x\in\liminf_n S_n$;
that is, $S\subset\liminf_n S_n$.
Now suppose that $y\notin S$. By (\ref{esb})
there is some $z\in\mathcal R^y\setminus S$.
Then $S\in V^{z}$. Consequently, there is some $m\in{\mathcal D}$ such
that $S_n\in V^{z}$ for all $n\succeq m$.
Hence, $S_n\cap\mathcal R_z=\varnothing$ for all $n\succeq m$.
Since $y\in\mathcal R_z$ and $\mathcal R_z$ is topologically open,
this implies
that $y\notin\limsup_n S_n$. Thus, $S\supset\limsup_n S_n$.
Summing it up, we have $\liminf_n S_n\supset S\supset\limsup_n S_n$.
Since $\liminf_n S_n\subset\limsup_n S_n$, this proves~(2).

To prove (3), suppose that $S=\lim_n S_n$. Then we know from (1) that
\mbox{$S\in\Scal$}.
Suppose that
$S_{n_j}$ is a subnet of $S_n$ that converges to some $S'\in\Scal$.
By~(2), we have $S'=\lim_j S_{n_j}$. Therefore, $S'=S$.
Thus, every convergent subnet of $S_n$ converges
to $S$. Since $\Scal$ is compact, this proves that $S_n$ converges to $S$,
and completes the proof of the lemma.\vspace*{2pt}
\end{pf}
\begin{lemma}\label{l2ndcountable}
Assuming (\ref{eopen}) and (\ref{esb}), if $X$ is second-countable
(i.e., with a countable base),
then $\Scal$ is second-countable and hence metrizable.
\end{lemma}
\begin{pf}
Assume that $X$ is second-countable.
In particular, there is\break a~countable dense
set $X_0\subset X$. Let $\mathcal B$ be a countable basis for $\topo_X$.
We claim that $\mathcal B':=\{V_U\dvtx U\in\mathcal B\}\cup\{V^x\dvtx
x\in X_0\}$
is a subbase for $\topo$.

Indeed,\vspace*{1pt} if $U'\in\topo_X$, then there is a subset $J\subset\mathcal
B$ such
that $U'=\bigcup_{U\in J} U$. Then $V_{U'}=\bigcup_{U\in J} V_U$; thus,
$V_{U'}$ is in the minimal topology containing $\mathcal B'$.

Now suppose that $x\in X$, then
%
%
\begin{equation}\label{eqvx}
V^x=\bigcup_{y\in X_0\cap\mathcal R^x} V^y.
\end{equation}
Indeed, if $S\in V^x$, then $x\notin S$ and
by (\ref{esb}) the topologically open set $\mathcal R^x\setminus S$ is
nonempty.
Take some $y\in X_0\cap\mathcal R^x\setminus S$, then $S\in V^y$.
This proves that $V^x\subset\bigcup_{y\in X_0\cap\mathcal R^x} V^y$.
Since\vspace*{1pt} $y<x$, we also have $V^y\subset V^x$ for $y\in\mathcal R^x$,
concluding that $V^x=\bigcup_{y\in X_0\cap\mathcal R^x} V^y$.
Thus, $V^x$ is also in the minimal topology containing $\mathcal
B'$.\vspace*{1pt}

Now Remark \ref{raltsubbase} shows that $\mathcal B'$ is a subbase
for $\topo$. Since $\mathcal B'$ is countable,~$\topo$
also has a countable basis (of all finite intersections of
elements of $\mathcal B'$).

By Proposition \ref{pHau}, topology $\topo$ is Hausdorff.
Since a second-countable compact Hausdorff space is metrizable
(which follows from Urysohn's Theorem 4.16 in \cite{kelley-book}),
this completes the proof.\vspace*{-2pt}
\end{pf}
\begin{lemma}\label{lsigma}
Suppose that (\ref{eopen}) and (\ref{esb}) hold,
and that $X_0$ is a dense subset of $X$.
Then the Borel $\sigma$-field of $(\Scal,\topo)$
is generated by $V^y, y\in X_0$.\vspace*{-2pt}
\end{lemma}
\begin{pf}
Observe that the Borel $\sigma$-field of $(\Scal,\topo)$
is generated by topologically open sets, which in turn are
generated by the sets of the form~$V_U$ and $V^x$
for topologically open $U\subset X$ and $x\in X$.
So we must prove that those can be obtained from $V^y, y\in X_0$, by
the operations of countable union, intersection and taking complements.

For $x\in X$, the set $V^x$ can be obtained as
$V^x=\bigcup_{y\in X_0\cap\mathcal R^x} V^y$,
by (\ref{eqvx}).

For a topologically open $U\subset X$, note that by the definitions
$V_U=\break \bigcup_{x\in U}\neg V^x$.
Therefore, $V_U\supset\bigcup_{y\in U\cap X_0}\neg V^y$.
On the other hand, take $x\in U$. By (\ref{esb}) we can find
$y\in X_0$ inside the topologically open set $U\cap\mathcal R^x$, then
$y< x$ and $\neg V^x\subset\neg V^y$.
Thus $V_U\subset\bigcup_{y\in U\cap X_0}\neg V^y$.
We conclude that
$V_U= \bigcup_{y\in U\cap X_0}\neg V^y$,
and so the sets $V_U$ can also be obtained.\vspace*{-2pt}
\end{pf}

Finally, observe that a topologically open $X'\subset X$
inherits the properties~(\ref{eopen}) and (\ref{esb}) from $X$,
and the following continuity under inclusions holds.\vspace*{-2pt}
\begin{lemma}\label{lmonot}
Let $X'\subset X$ be topologically open, and for $S\in\Scal_X$
let $\Phi(S):= S\cap{X'}$. Then $\Phi$ is a continuous map from
$\Scal_X$ to $\Scal_{X'}$.\vspace*{-2pt}
\end{lemma}
\begin{pf}
We have to check that preimages of topologically open sets are also
topologically open.
This follows easily since a topologically open subset $U\subset X'$
is topologically open inside\vspace*{1pt} $X$ and $\Phi^{-1}(V^U)=V^U$,
while for $x\in X'$ we have $\Phi^{-1}(V^x)=V^x$.\vadjust{\eject}
\end{pf}

Summing it up, we arrive at the following.
\begin{theorem}\label{tscal}
Let $(X,\tau)$ be a partially ordered second-countable topological space,
such that the ordering satisfies (\ref{eopen}) and (\ref{esb}).
Then:
\begin{longlist}[(2)]
\item[(1)] $(\Scal_X,\topo_\Scal)$ is a compact metrizable
Hausdorff space.
%
\item[(2)] Let $X_0$ be a dense subset of $X$.
If $S_1,S_2\in\Scal_X$ satisfy $S_1\cap X_0=S_2\cap X_0$, then $S_1=S_2$.
Moreover, the $\sigma$-field generated by $V^y$, $y\in X_0$, is
the Borel $\sigma$-field of $(\Scal_X,\topo_\Scal)$.
\item[(3)] If $\{S_n\}$ is a sequence in $\Scal_X$ and $S\in\Scal
_X$, then
$S_n\to S$ in $\topo_\Scal$ is equivalent to $S=\limsup_n
S_n=\liminf_n S_n$.
%
\item[(4)] Let $X'\subset X$ be a topologically open subset, and for
$S\in\Scal_X$
let $\Phi(S):= S\cap{X'}$. Then $\Phi$ is a continuous map from
$\Scal_X$ to $\Scal_{X'}$.
\end{longlist}
\end{theorem}
\begin{pf}
Part 1 follows from Propositions~\ref{pcompact},~\ref{plower} and
Lemma~\ref{l2ndcountable};
part~(2) follows from Lemmas~\ref{ldensedetermines} and~\ref{lsigma};
part~(3) follows from Lemma~\ref{llim};
and part~(4) follows from Lemma~\ref{lmonot}.
\end{pf}

\section{Uniformity in the mesh size}

In this section, we prove a mesh-indepen\-dent version of the gluing theorem.
\begin{proposition}[(Mesh-independent gluing)]\label{pmain}
Let $D\subset\C$ be some domain, and let $Q_0\in\Quad_D$
be a piecewise smooth quad in $D$.
Suppose that $\alpha\subset\C$ is a finite union of finite length paths,
with finitely many double points.
Consider a collection of percolation models indexed
by a set $\{\eta\}$.
Assume that~$\Pdisc_{\eta}$ converges
to (a subsequential) scaling limit $\Plim$
as mesh $|\eta|\to0$.
Then for every $\eps>0$ there is a finite set of piecewise smooth quads
$\Quad_\eps\subset\Quad_{D\setminus\alpha}$
and a subset ${\mathcal{W}}_\eps\subset\Scal$, which is measurable
with respect
to the finite $\sigma$-field $\sigma(\cross Q\dvtx Q\in\Quad_\eps)$,
such that
%
%
\begin{equation}
\label{efinapprox}
\lim_{|\eta|\to0} \Pdisc_\eta({\mathcal{W}}_\eps\Delta\cross
{Q_0})<\eps.
\end{equation}
\end{proposition}
%
%
\begin{pf}
This is the most technically difficult part of our paper.
The rough plan is as follows:
we set up a thin strip around $\alpha$,
and cut it into ``bays,'' forming a neighborhood of $\alpha$,
and a part away from $\alpha$, bounded by ``beaches.''
Partition depends on the percolation configuration,
and has the properties that
(i) percolation in the neighborhood of $\alpha$ is decoupled from the rest,
and
(ii) the crossings outside this neighborhood can be encoded by a finite
information,
stable under small perturbations of the percolation configuration.
We use a complicated construction of the neighborhood for the sake of
the property (ii).
To effectively bound the needed crossing information away from $\alpha$
(whose amount grows as we pass to the scaling limit), we further
coarse-grain this procedure to a fixed small scale,
when it is described by a finite $\sigma$-field.
This will allow to write the required estimates.

We now start the proof by fixing a countable subsequence $\{\eta_j\}$
from our set, such that $|\eta_j|\to0$ (and $\Pdisc_{\eta_j}$ converges
to $\Pdisc_0$).
Recall that by Remark~\ref{remseqnet} in our setting it
is sufficient to work with sequences, rather than nets.

\subsection*{Setup for the neighborhood of $\alpha$}
Reasoning like in the proof of Theorem~\ref{tdiscrete}
(the first paragraph of Section \ref{secdiscretegluing}),
we can approximate the quad~$Q_0$,
by quads whose boundaries intersect $\alpha$ on a finite set.
Thus, we can assume that $\alpha$
has finitely many double points, intersects ${\partial}Q_0$ on a
finite set.
Prolonging some of the curves in $\alpha$
if necessary (to cut the non-simply-connected components into smaller pieces),
we can further assume that each component of $[Q_0]\setminus\alpha$
is simply connected.\vadjust{\goodbreak}

Fix small $\eps_0>0$. 
Define two random variables: one given by the indicator function
of the crossing event $Y_0:= 1_{\cross{Q_0}}$;
and another by the conditional probability of a crossing, given the
percolation configuration $s$-away from $\alpha$:
\[
Y_s:=\Pdisc_\eta(\cross{Q_0}\mid{\mathcal{F}}_s).
\]
Here ${\mathcal{F}}_s$, $s>0$,
denotes the $\sigma$-field generated by the restriction of
$\omega$ to the complement of the $s$-neighborhood of $\alpha$,
as in Theorem \ref{tdiscrete}. By that theorem,
for all sufficiently small $s>0$,
%
%
\begin{equation}
\label{egoods}
\sup_{|\eta|\in(0,s)} \Pdisc_\eta(\eps_0<Y_s<1-\eps_0)<\eps
_0 .
\end{equation}
Denote by $x_i$, $i=1,\ldots,n$,
the finitely many points where $\alpha$ intersects ${\partial}Q_0$,
as well as the endpoints and the double points
of curves in $\alpha$.
Let $d$ be the minimum over $i$ of the distances from $x_i$ to other $x_j$'s
and the farthest of~${\partial}_0 Q_0$ and ${\partial}_2 Q_0$.
We now fix $s>0$ sufficiently small so that
(\ref{egoods}) holds, and sufficiently small
so that for $i=1,\ldots,n$, the balls $B(x_i,2s)$ are disjoint, and
%
%
\begin{equation}
\label{esonearm}
\sup_{|\eta|\in(0,s)}
\Pdisc_\eta\bigl(\exists\mbox{ a crossing from }B(x_i,2 s)\mbox{ to
}{\partial}B(x_i,d/3)\bigr)<\eps_0/n .
\end{equation}
The latter follows for small $s$ from the RSW estimate (\ref{eqwrsw}).
We also require that the intersection of $\alpha$ with $\bigcup
_i{\partial}B(x_i,s)$
is finite (which holds for almost every value of $s$,
since the total length of $\alpha$ is finite).

For the following construction, the reader is advised to refer to
Figure~\ref{figpaths}.
%
%
\begin{figure}

\includegraphics{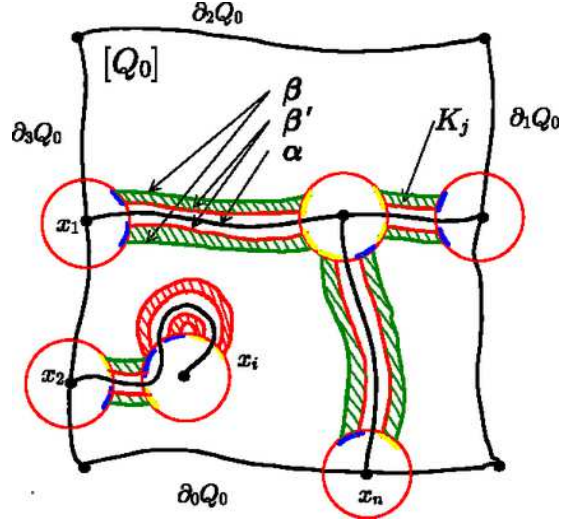}

\caption{Curves in $\alpha$ are between $\beta'$ which are in turn
between $\beta$.
We modify configuration on the small discs so that long quads bounded
by $\beta$ and $\beta'$
on long sides, have open tiles on one short side and closed on another.}
\label{figpaths}
\end{figure}
Let $\beta$ and $\beta'$ be finite\vadjust{\goodbreak} unions of disjoint smooth simple
paths in $[Q_0]\setminus\bigcup_iB(x_i,s)$,
and let $K$ be the closure of the union of connected components of
\[
[Q_0]\Bigm\backslash\biggl(\bigcup_i
B(x_i,s)\cup\beta\cup\beta'\biggr),
\]
whose boundary intersects both $\beta$ and $\beta'$.
We can choose these unions of paths so that:
\begin{itemize}
\item all their endpoints belong to $\bigcup_i {\partial}B(x_i,s)$,
\item$\beta'$ separates $\beta$ from $\alpha$,
\item the connected component of $[Q_0]\setminus\beta$ that contains
$\alpha$ is contained in the $(s/2)$-neighborhood of $\alpha$,
\item each component of connectivity $K_j$ of $K$
is a quad with two ``long'' sides
on $\beta$ and $\beta'$ and two ``short'' sides on two of the circles
${\partial}B(x_i,s)$.
\end{itemize}
It is immediate to verify that there exist such unions of paths.




\subsection*{Bays and beaches}
The following construction is illustrated in Figure \ref{figbays}.
Let $\omega$ be a sample from $\Pdisc_\eta$ for some small $\eta>0$.
For somewhat technical reasons, it will be convenient to consider in place
of $\omega$ the configuration $\tilde\omega$ which is modified on
the tiles intersecting disks $B(x_j,s)$.
We alter $\omega$ so that
for every component $K_i$ one of two short sides is covered by open tiles
and another by closed.
Clearly, (\ref{esonearm}) implies that if $\eta<s$,
%
%
\begin{equation}
\label{etildeom}
\Pdisc_\eta(\{\omega\in\cross{Q_0}\}\Delta\{\tilde\omega\in
\cross{Q_0}\})
\le\eps_0 .
\end{equation}
For this reason, studying $\tilde\omega$ will still be useful.

%
\begin{figure}

\includegraphics{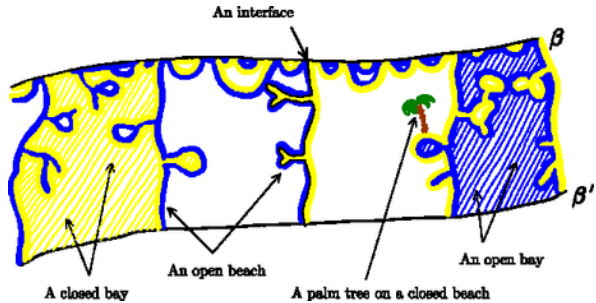}

\caption{Interfaces between open and closed clusters,
starting from $\beta$,
cut the strips $K_i$ 
into components.
\textit{Bays} are the components touching $\beta'$, and are
alternatively open and closed.
\textit{Beaches} are the inner boundaries of the bays, when entered
through $\beta'$.}
\label{figbays}
\vspace*{6pt}
\end{figure}

In what follows, we start on the curve $\beta$ and explore
all the potential crossings of the strip between $\beta$ and $\beta'$.
To be more precise,
denote by ${\partial}\tilde\omega$ the set of percolation
interfaces--curves, separating open tiles in
$\tilde\omega$ from the closed ones.
Let $I$ be the set of connected components of the intersection of~${\partial}\tilde\omega$
with the interior of $K$ (i.e., percolation interface in the interior
of $K$),
and let $\tilde I$ denote the set of elements of $I$ that have at least
one endpoint
on $\beta$.
The definition of $\tilde\omega$ guarantees that in each component $K_i$
there is at least one interface of $\tilde\omega$ with one endpoint
on $\beta$
and the other on $\beta'$.

Let $\Gamma$ denote the union of the elements of $\tilde I$; that is,
$\Gamma=\bigcup\tilde I=\bigcup_{\gamma\in\tilde I}\gamma$.
Each connected component of $K\setminus\Gamma$ which meets $\beta'$
will be called a
\textit{bay}.
Clearly, each point in $\beta'$ is in the closure of some bay,
and for each bay its intersection with $\beta'$ is connected.
Let $\bay$ be some bay.
The union of the tiles of $H_\eta$ that intersect $\bay$ as well as
${\partial}\bay\setminus\beta'$ will be called the \textit{beach} of
$\bay$.
It is easy to see that the beach is connected, and all the tiles
in the beach are either open (i.e., in $\tilde\omega$) or closed (not
in~$\tilde\omega$).
If they are all open (resp., closed), then $\bay$ itself will be
referred to as an open (resp., closed) bay.

Let $M'=M'(\omega)$ be the union
of the component of $[Q_0]\setminus\beta'$ containing $\alpha$,
and of all the bays without beaches.
Denote by $M=M(\omega)$ the complement in $[Q_0]$ of $M'$,
as well as the corresponding set of tiles.
The following property of $M$ is essential: given $M$
and the restriction of $\omega$ to it, 
the conditional law of the restriction
of $\omega$ to $M'$ (i.e., the complement of $M$) is unbiased.
In other words, if we take $\omega_1$ independent from $\omega$
and of the same law, and we define $\omega_2$ to agree with
$\omega$ on $M$ and to agree with
$\omega_1$ elsewhere, then $\omega_2$ also has the law $\Pdisc_\eta$.
This follows directly from the fact that $M(\omega_2)=M(\omega)$.
Below we will denote the restriction of a percolation configuration
$\omega$ to $M$ by $\omega\llcorner M$.

Fix some small $s'\in(0,s)$,
and let $\exep=\exep_{s'}$ denote the event that there is a~small bay
$\bay$,
that is, with $\diam(\bay\cap\beta')<s'$.
Observe that for every bay $\bay$ the endpoints of $\bay\cap\beta'$
are endpoints of interfaces of $\omega$ which meet both $\beta$ and
$\beta'$.
Thus, three crossings of $K$ (which are alternating, e.g., closed, open
and closed) land on the arc $\bay\cap\beta'$
and Lemma \ref{lthree} applied to all components of $K$ implies that
\[
\lim_{s'\to0}
\sup_{|\eta|\in(0,s')} \Pdisc_\eta(\exep) = 0.
\]
By choosing $s'$ sufficiently small, we therefore assume, with no loss
of generality, that $\Pdisc_\eta(\exep)<\eps_0$ holds for
every $\eta\in(0,s')$ and that $s'<s$.
Observe that on $\neg\exep$ the number of bays is bounded by a constant
depending only on~$s',\beta,\beta'$ and $Q_0$.

We will find the structure of the bays to be a convenient component
in the estimation of the conditional probability (given the
bays and some additional information) of $\tilde\omega\in\cross{Q_0}$.
However, as $|\eta|\to0$, the number of possible bay
configurations increases without bound,
and therefore their structure cannot be completely captured by a
finite $\sigma$-field. For this purpose, we need to introduce a discretization,
and any sort of coarse-graining, for example, by the square lattice,
would do.

\vspace*{5pt}\subsection*{Discretization to a finite $\sigma$-field by coarse-graining}

We need some tessellation of $K$, and
the following one is somewhat ad hoc, but suitable.
Let $\delta\in(0,s')$ be small, to be chosen later.
Fix some finite set of points $W\subset K$
such that the distance between any two points in $W$ is
at least $\delta/4$ and every point in $K$ is within distance $\delta/2$
of some point in $W$. Let $T$ be the collection of Voronoi tiles
in $K$ with respect to $W$; that is, each $T_w$, $w\in W$, is
the set of points $x\in K$ whose distance to $W$ is equal to $|x-w|$.
Then the diameters of the tiles in the tessellation $T$ are all less
than $\delta$
(but the lattice~$H_\eta$ tiles are still much smaller).
Let $\hat T$ denote the tessellation obtained by adjoining to $T$
the closures of the (finite in number) connected components of
$[Q_0]\setminus K$ that do not contain $\alpha$.
We think of tessellations $T$ and $\hat T$ as cell complexes.
The vertices of $T$ are defined as the vertices of the Voronoi tiles
together with
the points in $(\beta\cup\beta')\cap{\partial}{Q_0}$.
The vertices of~$\hat T$ are defined as those of~$T$ plus the four
corners of $Q_0$.
Note that some of the edges of $\hat T$ are not necessarily straight.
Denote by $[\hat T]$ the union of the tiles in $\hat T$.

We say that a quad $Q\in\Quad_D$ is \textit{compatible} with $\hat T$ if
$[Q]$ is a union of tiles of $\hat T$ and
each of its four corners $Q(0,0)$, $Q(0,1)$, $Q(1,0)$, $Q(1,1)$
is a vertex of~$\hat T$.
Since $\hat T$ is finite, it is clear that up to reparameterization
preserving corners there are only finitely
many compatible $Q$. Let $\Quad_{\hat T}$ denote a set of
representatives of equivalence
classes of compatible quads up to reparameterization preserving corners.
Then $\Quad_{\hat T}$ is a finite collection of
piecewise smooth quads in $D\setminus\alpha$.
Let ${\mathcal{F}}_T$ denote the $\sigma$-field generated by the events
$\cross Q$, $Q\in\Quad_{\hat T}$.

Without loss of generality, we assume
that along our chosen sequence~$\{\eta_j\}$, for
each edge $e$ of the tessellation $T$ that does not touch
${\partial}{Q_0}$:
\begin{itemize}
\item each edge of $H_\eta$ intersects the tessellation edge $e$ in at
most finitely many points,
\item the tessellation edge $e$ does not contain any vertices of
$H_\eta$,
\item the endpoints of the tessellation edge $e$ are not on edges of
$H_\eta$.
\end{itemize}
There is no loss of generality in this assumption, because
we may slightly perturb $\beta$, $\beta'$ and $T$
(and we do not even need $T$ to be a Voronoi tessellation).

Let $\bay$ be some bay. Define $T(\bay)$ to be its $T$-discretization,
namely the union of all edges of $T$ inside $\bay\cap\beta'$
and all tiles of $T$ inside $\bay$, which can be connected to
$\beta'$ by a path of tiles of $T$ also contained inside $\bay$.
By construction, each point in ${\partial}T(\bay)$ is within distance
$\delta$ from
${\partial}\bay$.
Note that there is a finite number of possibilities
for the choices of $T(\bay)$,
and also observe
that whenever $\exep$ does not hold,
all $T$-discretizations are nonempty and so
$T(\bay)=T(\bay')$ implies $\bay=\bay'$.

\subsection*{Connectivity in the coarse-grained model}
The next objective is to show that on the complement of $\exep$ the
collections of discretized bays
\[
Z_o := \{T(\bay)\dvtx\bay\mbox{ an open bay}\}
\]
and
\[
Z_c := \{T(\bay)\dvtx\bay\mbox{ a closed bay}\}
\]
can be determined from ${\mathcal{F}}_T$.
That is, there are ${\mathcal{F}}_T$-measurable
random variables $Z_o'$ and $Z_c'$ such that $Z_o=Z'_o$
and $Z_c=Z'_c$ on the complement of $\exep$.
Here and in the following, we ignore events that have zero
$\Pdisc_\eta$-measure for each $\eta\in\{\eta_j\}$.

Note that on the event $\neg\exep$ the cardinality of $Z_o$ and $Z_c$
is bounded by a~constant depending on $Q_0$, $\beta'$ and $s'$ only.

Let $e$ be an edge of $T$ that lies on the component $\beta'_j=\beta
'\cap K_j$,
\mbox{denote} also $\beta_j:=\beta\cap K_j$.
Let $Q_e$ and $Q'_e$ be the quads in $\Quad_{\hat T}$ that satisfy
${\partial}_0 Q_e=\allowbreak{\partial}_1 Q'_e=e$,
${\partial}_2 Q_e={\partial}_3 Q'_e=\beta_j$ and $[Q_e]=[Q'_e]=K_j$.
Then the event that there is an interface $\gamma\in\tilde I$
that meets $e$ is the same as the event $\cross{Q_e}\setminus\cross{Q'_e}$.

If $E$ is the set of all edges of $T$ on $\beta'_j$ which meet
interfaces in $\tilde I$,
then on the event $\neg\exep$,
the connected components of $\beta'_j\setminus\bigcup E$ are the
intersections of the elements of $Z_o\cup Z_c$ with $\beta'_j$,
and since the bays alternate in color, it follows that
the sets
$\{\tau\cap\beta_j'\dvtx\tau\in Z_o\}$
and
$\{\tau\cap\beta_j'\dvtx\tau\in Z_c\}$
can each be determined from ${\mathcal{F}}_T$, up to an event
contained in
$\exep$.

Next, suppose that $\gamma\subset\beta_j'$ is of the form
$\tau_0\cap\beta_j'$ for some $\tau_0\in Z_c$
(assuming~$\neg\exep$, such intersection is always nonempty).
Let $T_0$ be some tile of $T$.
Then $T_0\subset\tau_0$ if and only if there is a simple path $\gamma'$
of edges of $T$ with the endpoints of $\gamma'$ on $\gamma$
such that $\gamma'$ separates
$T_0$ from $\beta_j$ in $K_j$ and there is no open crossing
in $K_j\cap\tilde\omega$ from $\gamma'$ to $\beta_j$.
This, along with the dual argument,
implies the existence of $Z_o'$ and $Z_c'$, as claimed.

\subsection*{Describing the crossing structure by a finite graph}
We now introduce a~graph that describes the connectivity
(by open crossings)
structure of the various open beaches and the two boundary edges
${\partial}_0Q_0$ and ${\partial}_2 Q_0$.

We start by defining a random graph $G$
which would describe the connectivity away from $\alpha$.
The vertices of $G$ are
$Z_o\cup\{{\partial}_0Q_0,{\partial}_2Q_0\}$, where an edge is placed
between vertices $v,v'\in Z_o$ if the beaches of the corresponding
bays are connected by an open crossing in
the restriction of $\tilde\omega$ to $M$, denoted by
$\tilde\omega\llcorner M$;
while if $v$ and/or $v'$ are in $\{{\partial}_0Q_0,{\partial}_2Q_0\}
$, the connectivity to the beach
is simply replaced by the connectivity to the corresponding boundary edge
itself.
Note that there is a finite number of possibilities for the choice of $G$.

Suppose that $\tau_0,\tau_1\in Z_o$. Let $e_j$, $j=0,1$, be an edge
of $T$ that is contained in $\beta'$ and has an endpoint in $\tau_j$,
but $e_j$ itself is not in $\tau_j$.
Then $e_j$ meets an interface
on the boundary of the corresponding bay.
On $\neg\exep$, the two edges $e_j$ are connected
by an open crossing in $\tilde\omega\llcorner[\hat T]$ if
and only if $[\tau_0,\tau_1]$ is an edge of $G$.
Thus, on the event $\neg\exep$ the subgraph of $G$ induced by $Z_o$
is determined from ${\mathcal{F}}_T$. A similar argument applies also
to the
edges with endpoints in $\{{\partial}_0Q_0,{\partial}_2Q_0\}$. Thus,
on the complement of $\exep$ the graph $G$ can be determined from
${\mathcal{F}}_T$.

Now we define a random graph $G^*$ so that it describes the
connectivity near~$\alpha$.
Consider $L=({\partial}_0Q_0\cup{\partial}_2Q_0)\setminus[\hat T]$,
which contains finitely many arcs on ${\partial}Q_0$.
Denote by $G^*$ the graph whose vertices are $Z_o$ and the connected
components of $L$, where $v,v'\in Z_o$ are connected by an edge in $G^*$
if the two corresponding beaches of $\tilde\omega$ are connected
by an open crossing in
\mbox{$\tilde\omega\llcorner([Q_0]\setminus M)$}, and if $v$ and/or
$v'$ are
components of $L$, then the connectivity to the beaches is replaced by
connectivity to $v$ and/or $v'$.

Clearly, $\tilde\omega\in\cross{Q_0}$ if and only if there is
a path from ${\partial}_0Q_0$ (or a subarc thereof) to ${\partial
}_2Q_2$ in
$G\cup G^*$.
Thus, on the event $\neg\exep$,
the graphs $G$ and $G^*$ determine whether $\tilde\omega\in\cross{Q_0}$.

Denoting by $\omega\llcorner M$ the restriction of $\omega$ to $M$,
we set
\[
\tilde Y=\tilde Y(\omega):= \Pdisc_\eta(\tilde\omega\in\cross
{Q_0}\mid\omega\llcorner{M}) .
\]
We now show that on $\neg\exep$, the knowledge of $G$, $Z_c$ and $Z_o$
can be used to approximate $\tilde Y$.

Let
\[
T(M):=[\hat T]\bigm\backslash\Bigl(\bigcup Z_c\cup\bigcup Z_o\Bigr) ,
\]
that is, the complement in $[\hat T]$ of the union
of $T(\bay)$, where $\bay$ runs over all the bays.

Let $G^\dagger$ be the graph on the same vertex set as $G^*$, where
an edge appears between $v,v'\in Z_o$ if the common boundary of $v$
and $T(M)$ is connected in $\tilde\omega\llcorner([Q_0]\setminus
T(M))$ with the
common boundary of $v'$ and $T(M)$,
while if~$v$ and/or $v'$ are components of $L$, then the connectivity
is instead to $v$ and/or~$v'$, but still within
$\tilde\omega\llcorner([Q_0]\setminus T(M))$.
We now show that if $\delta$ is sufficiently small, $\eta<\delta$
and $\eta\in\{\eta_n\}$,
then
%
%
\begin{equation}
\label{eGG}
\Pdisc_{\eta}(G^\dagger\ne G^*, \neg\exep\mid\omega\llcorner
M)<\eps_0 .
\end{equation}
As we have mentioned, on the event $\neg\exep$ there is a finite
upper bound on the size of $Z_c\cup Z_o$,
depending on $Q_0$, $\beta'$ and $s'$ only.
Therefore, it is enough to have a good estimate for the probability that
a particular pair of vertices $v$ and $v'$ of $G^*$ are connected by an
edge in one
of the two graphs but not in the other,
and we are free to impose additional conditions on $\delta$.
Connectivity of any pair of vertices
is reduced to $\delta$-perturbations of quads of size at least~$s'$,
and so is easily obtained by several applications
of Lemma \ref{lcont} (some of them applied with the corners of the
quads appropriately
permuted) proving~(\ref{eGG}).

\subsection*{Final estimates}
Define $\tilde Y_0=\tilde Y_0(\omega):=1_{\tilde\omega\in\cross
{Q_0}}$ and
\[
\tilde Y_T:=
\Pdisc_{\eta}({\partial}_0 Q_0 \mbox{ and }{\partial}_2 Q_0\mbox{
are connected
by a path in }G^\dagger\cup G\mid\omega\llcorner M).
\]
Then $\tilde Y_T$ is a function of the triple $(Z_c,Z_o,G)$,
and is therefore ${\mathcal{F}}_T$-measura\-ble.
By the above observation that
on $\neg\exep$ the event $\cross{Q_0}$ is described by connectivity
in $G\cup G^*$, we have
\[
|\tilde Y_T-\tilde Y| 1_{\neg\exep}
\le
\Pdisc_{\eta}(G^\dagger\ne G^*, \neg\exep\mid\omega\llcorner M).
\]
By taking expectations, applying (\ref{eGG}) and recalling that
$s'$ was chosen to guarantee $\Pdisc(\exep)<\eps_0$ one gets
%
%
\begin{equation}
\label{eYY}
\|\tilde Y_T-\tilde Y\|_2<\sqrt{2\eps_0} ,
\end{equation}
where the norm refers to the measure $\Pdisc_\eta$.

By (\ref{egoods}) and the definition of $Y_s$, we can write
\[
\|Y_0-Y_s\|_2^2<2 \eps_0+\eps_0^2<3 \eps_0 ,
\]
and by (\ref{etildeom}), we have
%
%
\begin{equation}
\label{eyy2}
\|Y_0-\tilde Y_0\|_2<\sqrt{\eps_0} . 
\end{equation}
Hence, $\|\tilde Y_0-Y_s\|_2<\sqrt{3 \eps_0}+\sqrt{\eps_0}<3 \sqrt
{\eps_0}$.
Since $Y_s$ is ${\mathcal{F}}_s$-measurable, it is also $\omega
\llcorner M$-measurable.
By its definition, $\tilde Y$ minimizes
$\| \tilde Y_0-X\|_2^2$ among $\omega\llcorner M$-measurable random
variables $X$.
Therefore, comparison to $Y_s$
(recall that the set $M$ includes the complement
of the $s$-neighborhood of $\alpha$)
yields
\[
\|\tilde Y_0-\tilde Y\|_2\le\|\tilde Y_0- Y_s\|_2<3 \sqrt{\eps_0} .
\]
Combining this with (\ref{eYY}) and (\ref{eyy2}),
we conclude that
\[
\|Y_0-\tilde Y_T\|_2 < 6 \sqrt{\eps_0} .
\]
Because $\tilde Y_T$ is ${\mathcal{F}}_T$-measurable
and $\eps_0$ may be chosen arbitrarily small, this
proves that there is an ${\mathcal{F}}_T$-measurable event
${\mathcal{W}}$ (which may depend on~$\eta$),
such that $\Pdisc_\eta({\mathcal{W}}\Delta\cross{Q_0})<\eps/2$.
However, since ${\mathcal{F}}_T$ is finite,
there are finitely many possibilities
for the event ${\mathcal{W}}$, and one of those
works for a sub-subsequence $\eta_{j_k}$.
Since we work with a (subsequential) scaling limit,
the limit in (\ref{efinapprox}) exists
along the original subsequence and is equal to $\Pdisc_0({\mathcal
{W}}\Delta\cross{Q_0})$.
Since along a sub-subsequence the quantity in question is bounded by
$\epsilon$,
so is the limit and we deduce (\ref{efinapprox}).
\end{pf}

\section{Factorization}

$\!\!\!$In this section, we prove the Factorization
Theorem~\ref{tfactor}.\vspace*{2pt}
\begin{lemma}\label{lcrossingOpen}
Let $\Pdisc_0$ be some subsequential scaling limit, then
the boundary (in topology $\topo$ on $\Scal$) of a crossing event has
probability zero.
Namely,
\[
\Plim({\partial}\cross{Q_0})=0
\]
holds for every $Q_0\in\Quad_D$.
\end{lemma}
\begin{pf}
Fix $Q_0\in\Quad_D$ and let $\eps>0$.
It is easy to see (e.g., using the Riemann map onto $\hat\C\!\setminus[Q_0]$)
that there is a continuous injective\vspace*{2pt} \mbox{$\hat Q_0\dvtx[-1,2]^2\,{\to}\,\C$}, whose
restriction
to $[0,1]^2$ is $Q_0$.
Let $M:=[-1,1/3]^2\times[2/3,2]^2$.
For every quadruple $q=(x_0,y_0,x_1,y_1)\in M$, define the quad
$Q^q\dvtx[0,1]^2\to\C$ by
\[
Q^q(x,y):=\hat Q_0\bigl(x_0+(x_1-x_0) x,y_0+(y_1-y_0) y\bigr).
\]
It is\vspace*{2pt} a perturbation of $Q_0$, obtained as an image of $[x_0,x_1]\times
[y_0,y_1]$ by $\hat Q_0$. Lemma~\ref{lcont} implies that there is some positive $\delta_0>0$
(depending on $\hat Q_0$), such that if $q$ and $q'$ are in $M$ and differ in exactly one
coordinate, and the difference in that coordinate is at most $\delta_0$,
then $\limsup_{|\eta|\to0}\Pdisc_\eta(\cross{Q^q}\Delta\cross{Q^{q'}})<\eps$.
Consequently, we have
\[
\limsup_{|\eta|\to0}\Pdisc_\eta(\cross{Q^q}\Delta\cross
{Q^{q'}})<4 \eps\qquad
\mbox{provided }\|q-q'\|_\infty\le\delta_0.
\]
Set $q_s:=(-s,s,1+s,1-s)$, $Q':=Q^{q_{-\delta_0/2}}$
and $Q'':=Q^{q_{\delta_0/2}}$.
Then $Q'<Q_0<Q''$ and
%
%
\begin{equation}
\label{eQQdiff}
\limsup_{|\eta|\to0}\Pdisc_\eta(\cross{Q'}\Delta\cross
{Q''})<4 \eps.
\end{equation}
Also, $Q',Q''\in\Quad_D$, if $\delta_0$ is chosen sufficiently small.

Recall that $V^{Q_0}=\Scal_D\setminus\cross{Q_0}$ is open in $\Scal_D$.
Hence $\cross{Q_0}$ is closed.
Let $U'$ be the set of quads $Q\in\Quad_D$ satisfying $Q_0<Q<Q''$.
Then
\[
\cross{Q_0}\supset V_{U'}\supset\cross{Q''}.
\]
By passing to the complements, we see that there is a
closed subset (the complement of $V_{U'}$) that contains
$V^{Q_0}=\neg\cross{Q_0}$ and is contained in $V^{Q''}$.
Hence, the closure of $V^{Q_0}$ is contained in $V^{Q''}$, which gives
\[
{\partial}\cross{Q_0}={\partial}V^{Q_0} \subset\closure{V}{}^{Q_0}
\subset V^{Q''}.
\]
On the other hand, let $U$ be the set of quads $Q\in\Quad_D$
satisfying $Q'<Q<Q''$, then
\[
{\partial}\cross{Q_0}\subset\cross{Q_0}\subset V_U.
\]

Now observe that $V_U\subset\cross{Q'}$ while $V^{Q''}= \neg\cross{Q''}$.
We have shown that the open set $V_U\cap V^{Q''}$ contains ${\partial
}\cross{Q_0}$.
The portmanteau theorem (see, e.g., \cite{dudley-book}, Theorem~11.1.1)
therefore gives
\begin{eqnarray*}
\Plim({\partial}\cross{Q_0})&\le&\Plim(V_U\cap V^{Q''})
\le\liminf_{|\eta|\to0}\Pdisc_\eta(V_U\cap V^{Q''})\\
&\le&\liminf_{|\eta|\to0}\Pdisc_\eta(\cross{Q'}\Delta\cross
{Q''}) .
\end{eqnarray*}
Since $\eps$ was arbitrary, the result now follows from (\ref{eQQdiff})
\end{pf}
\begin{corollary}\label{cpasslimit}
If ${\mathcal{M}}\subset\Scal_D$ is in the Boolean algebra generated by
finitely many of the
events $\cross Q$, $Q\in\Quad_D$
(i.e., it can be expressed using finitely many $\cross Q$ and the
operations of
union, intersection and taking complements),
then $\Plim({\mathcal{M}})=\lim_{|\eta|\to0}\Pdisc_\eta
({\mathcal{M}})$.
\end{corollary}
\begin{pf}
Lemma \ref{lcrossingOpen} implies that $\Plim({\partial}{\mathcal{M}})=0$.
Thus, the corollary follows from the portmanteau theorem \cite{dudley-book},
Theorem 11.1.1, and
the weak convergence of $\Pdisc_\eta$ to $\Plim$.
\end{pf}
\begin{pf*}{Proof of Theorem \ref{tfactor}}
Clearly, ${\mathcal{F}}_{D\setminus\alpha}=\bigvee_j {\mathcal{F}}_{D_j}$.
It therefore remains to prove the left-hand equality stated in the theorem.
Note that we work up to sets of $\Plim$-measure zero.

Take a smooth (given by a diffeomorphism) quad $Q_0\in\Quad_D$.
By Proposition \ref{pmain} and Corollary \ref{cpasslimit},
we have
%
%
\begin{equation}
\label{eQ0fac}
\cross{Q_0}\in{\mathcal{F}}_{D\setminus\alpha} .
\end{equation}
%
Since such collection of quads is dense in $\Quad_D$,
Theorem \ref{tfactor} follows from Proposition
\ref{tpscal}(2).
\end{pf*}

\begin{appendix}
\section{Continuity of crossing events}

In this section, we apply Russo--Seymour--Welsh techniques
to deduce estimates of various crossing probabilities.
We start by showing that crossing events
are stable under small perturbations of quads.
\begin{lemma}\label{lcont}
There exist a positive function $\rsww_c(\delta,d)$,
such that
\[
\lim_{\delta\to0}\rsww_c(\delta,d)=0 \qquad\mbox{for any fixed } d ,
\]
and the following estimates hold.
Let $Q$ be a quad.
Let $d_j$, $j=0,1$, be the infimum diameter of any path in $[Q]$ connecting
${\partial}_j Q$ and ${\partial}_{j+2}Q$, and define $d:=\max\{
d_0,d_1\}$.
Fix some $\delta<d/2$.
Let $Q'$ be another quad, satisfying at least
one of the following conditions;
%
%
%
\begin{figure}

\includegraphics{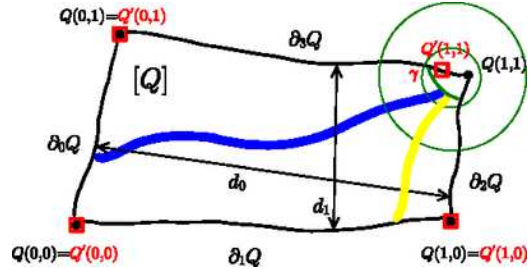}

\caption{Case \protect\hyperlink{casei}{(1)}: quads $Q$ and $Q'$ differ by only
one vertex,
separated from the other vertices by a short cut $\gamma$.
Configurations with only one of two quads crossed have
an open crossing from $\gamma$ to ${\partial}_0Q$ and a dual closed
crossing from $\gamma$ to $p_1Q$,
making their probability small.}
\label{figapprox1}
\vspace*{6pt}
\end{figure}
see Figures \ref{figapprox1} and \ref{figapprox2} for a graphical interpretation.
\begin{longlist}[(2)]
\item[(1)]
\hypertarget{casei}
$[Q']=[Q]$, ${\partial}_0Q'={\partial}_0Q$,
${\partial}_1Q'={\partial}_1 Q$, and there is
%
%
\begin{figure}[b]

\includegraphics{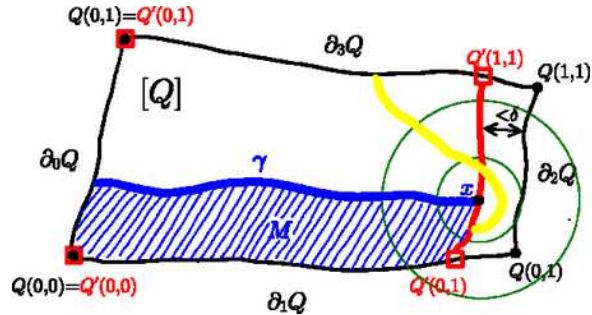}

\caption{Case \protect\hyperlink{caseii}{(2)}: side ${\partial}_2 Q'$ is close
to the side ${\partial}_2 Q$.
We take the lowest crossing $\gamma$ of $Q'$, landing at $x$.
A dual closed crossing prevents it from landing on ${\partial}_2Q$, making
the probability small.
Case \protect\hyperlink{caseiii}{(3)} is symmetric.}
\label{figapprox2}
\end{figure}
a path $\gamma\subset[Q]$ of diameter at most $\delta$ that separates
$\{Q(1,1),Q'(1,1)\}$ from ${\partial}_0Q\cup{\partial}_1 Q$ inside~$[Q]$.
\item[(2)]\hypertarget{caseii} $[Q']\subset[Q]$, ${\partial}_0Q'={\partial
}_0Q$, ${\partial}_1Q'\subset{\partial}_1 Q$,
${\partial}_3Q'\subset{\partial}_3 Q$, and each point on~${\partial
}_2Q'$ can be connected to ${\partial}_2 Q$
by a path $\alpha\subset[Q]$ with $\diam(\alpha)\le\delta$.
\item[(3)]\hypertarget{caseiii} $[Q']\subset[Q]$, ${\partial}_0Q'\subset
{\partial}_0Q$, ${\partial}_1Q'={\partial}_1 Q$,
${\partial}_2Q'\subset{\partial}_2 Q$, and each point on~${\partial
}_3Q'$ can be connected to ${\partial}_3 Q$
by a path $\alpha\subset[Q]$ with $\diam(\alpha)\le\delta$.
\end{longlist}
Then for every $|\eta|<\delta$ we have
\[
\label{emudiff}
\Pdisc_\eta(\cross{Q}\Delta\cross{Q'}) < \rsww_c(\delta,d) .
\]
\end{lemma}

For site percolation on the triangular lattice,
there is a short proof based on Cardy's formula.
The proof below is a simple application of the RSW estimate (\ref{eqwrsw})
and the ``lowest crossing'' concept.
\begin{pf*}{Proof of Lemma \ref{lcont}}
First, we deal with case \hyperlink{casei}{(1)} in the statement
and prove the estimate with $\rsww_c$ equal to $\rsww_1$
from assumption (\ref{eqwrsw}).
Suppose first that $d=d_0$.
The event $\cross Q\Delta\cross{Q'}$ is contained in the event
that there is a~percolation cluster meeting $\gamma$
and ${\partial}_0Q$.
The latter has an open crossing of the annulus $A(Q(1,1),\delta,d_0)$,
which by the RSW estimate (\ref{eqwrsw}) has probability at most
$\rsww_1(\delta,d)$,
and we are done with case \hyperlink{casei}{(1)}.
If $d=d_1$, a similar argument shows the symmetric difference between
the event
of a~closed crossing from~${\partial}_1Q$ to ${\partial}_3Q$ in $[Q]$
and the corresponding
event in $Q'$ is bounded by $\rsww_1(\delta,d_1)$.
Duality shows that the latter symmetric difference is
the same as $\cross Q\Delta\cross{Q'}$, which proves case~\hyperlink{casei}{(1)}.

We now deal with case \hyperlink{caseii}{(2)},
when clearly $\cross Q\Delta\cross{Q'}=\cross{Q'}\setminus\cross Q$.
We start by estimating the crossing probability for $Q'$.
Recall our convention that $\rsww(r,R)=1$ for $r\ge R$
(i.e., when annulus in question is empty).
Let $x$ be some point on the path $\tau$ of diameter $d_1$, connecting
${\partial}_1Q$ to ${\partial}_3Q$.
Any crossing of $Q'$ has diameter at least $d_0-\delta$ and passes
within distance
$\delta$ of $\tau$, so in particular it crosses an annulus
$A(x,d_1+\delta,(d_0-\delta)/2)$
if well defined,
and by the RSW estimate (\ref{eqwrsw}) we have
%
%
\begin{equation}\label{eqd1d}
\Pdisc_\eta(\cross{Q'}) \le\rsww_1\biggl(d_1+\delta,\frac{d_0-\delta
}{2}\biggr).
\end{equation}

Next, we cut ${\partial}_2 Q'$ by a point $z$ into two parts, $\sigma
_1$ and $\sigma_3$,
so that $\sigma_j$ cannot be connected to ${\partial}_j Q$ by a path of
diameter less than $d_1/2$ inside~$Q$.
If $\cross{Q'}$ occurs, there is a crossing to at least one of $\sigma
_1$ and $\sigma_3$.
We will work with configurations with a crossing landing on $\sigma_3$,
the other case being symmetric, with the lowest crossing replaced by
the uppermost,
so the total estimate will be double of what we obtain.

Consider a percolation configuration $\omega$.
On the event $\cross{ Q'}$,
let $\gamma$ be the ``lowest'' $\omega$-crossing of $Q'$ (in the sense
that it is the closest to ${\partial}_1 Q'$, i.e.,
separates ${\partial}_1 Q'$ from all other $\omega$-crossings of $Q'$ within
$[Q']$), and let $x$ be the endpoint of $\gamma$ on ${\partial}_2 Q'$.
Our assumption implies that $x\in\sigma_3$, or there would be a lower
crossing.

We now estimate $\Pdisc_\eta(\neg\cross Q\mid\cross{Q'},\gamma
,x\in\sigma_3)$.
Let $M$ be the connected component of $[Q']\setminus\gamma$ which
has ${\partial}_1 Q'$ on its boundary.
The event $\neg\cross Q$ would imply that $\gamma$ cannot
be connected to ${\partial}_2Q$ by a crossing inside $[Q]\setminus M$.
Hence, there is a dual closed crossing from ${\partial}_3 Q$ to
$(\sigma_3\cap\bar M)\cup{\partial}_1Q$,
in particular crossing the annulus $A(x,\delta,d_1/2)$
inside $[Q]\setminus M$.

The lowest crossing $\gamma$ depends only on the configuration inside $M$,
so the restriction of $\omega$ to $[Q]\setminus M$ is unbiased.
Therefore, the conditional probability of the mentioned annulus
crossing is
can be bound by the RSW estimate~(\ref{eqwrsw})
and we conclude that
%
%
\begin{equation}\label{eqdcond}
\Pdisc_\eta(\neg\cross Q\mid\cross{Q'},\gamma,x\in\sigma_3)
\le\rsww_1(\delta,d_1/2) .
\end{equation}

Now we are ready to prove the estimate, working out separately three
possibilities:
\begin{longlist}
\item $d=d_1$,
\item $d=d_0>\delta\ge d_1$,
\item $d=d_0>d_1>\delta$.
\end{longlist}
When (i) occurs, we majorate the probability
of the event in question by the maximum of its conditional probability,
estimated by (\ref{eqdcond}), to arrive at
%
%
\begin{equation}\label{eqconti}
\Pdisc_\eta(\cross{Q'}\setminus\cross Q)
\le2 \rsww_1\biggl(\delta,\frac{d_1}{2}\biggr)
= 2 \rsww_1\biggl(\delta,\frac{d}{2}\biggr) .
\end{equation}
When (ii) occurs, we use the estimate (\ref{eqd1d}):
%
%
\begin{equation}\label{eqcontii}
\Pdisc_\eta(\cross{Q'}\setminus\cross Q)
\le\Pdisc_\eta(\cross{Q'})
\le\rsww_1\biggl(d_1+\delta,\frac{d_0-\delta}{2}\biggr)
\le\rsww_1\biggl(2\delta,\frac{d}{4}\biggr) .
\end{equation}
Finally, when (iii) occurs, we apply the total expectation law,
using both estimates (\ref{eqdcond}) and (\ref{eqd1d})
along the way:
%
%
\begin{eqnarray}\label{eqcontiii}
\Pdisc_\eta(\cross{Q'}\setminus\cross Q)
&\le&\E[\Pdisc_\eta(\neg\cross Q\mid\cross{Q'},\gamma,x\in
\sigma_3)]\nonumber\\
&&{} +\E[\Pdisc_\eta(\neg\cross Q\mid\cross{Q'},\gamma,x\in
\sigma_1)]\nonumber\\
&\le& 2 \Pdisc_\eta(\cross{Q'}) \rsww_1\biggl(\delta,\frac
{d_1}{2}\biggr)\\
&\le& 2 \rsww_1\biggl(d_1+\delta,\frac{d_0-\delta}{2}\biggr) \rsww_1\biggl(\delta
,\frac{d_1}{2}\biggr)\nonumber\\
&\le& 2 \rsww_1\biggl(2 d_1,\frac{d}{4}\biggr) \rsww_1\biggl(\delta,\frac
{d_1}{2}\biggr).\nonumber
\end{eqnarray}
To prove lemma in case \hyperlink{caseii}{(2)}, we have to show that
for fixed $d$ and any $\epsilon>0$, the right-hand side of the estimates
(\ref{eqconti}), (\ref{eqcontii}) and (\ref{eqcontiii})
is smaller than~$\epsilon$, if $\delta$ is small enough.
For estimates (\ref{eqconti}) and (\ref{eqcontii}),
this follows directly from Assumption \ref{eqwrsw}.
For the remaining (\ref{eqcontiii}), let $\rho>0$ be such that
for $r<\rho$ we have $\rsww_1(2 r,\frac{d}{4})<\epsilon/2$.
Then if $d_1\le\rho$,
then the right-hand side of~(\ref{eqcontiii})
can be bounded by
\[
2 \rsww_1\biggl(2 d_1,\frac{d}{4}\biggr) \rsww_1\biggl(\delta,\frac{d_1}{2}\biggr)
\le2 \rsww_1\biggl(2 d_1,\frac{d}{4}\biggr) < \epsilon,
\]
we are done.
Otherwise, $d_1>\rho$ and the right-hand side of (\ref{eqcontiii})
can be bounded by
\[
2 \rsww_1\biggl(2 d_1,\frac{d}{4}\biggr) \rsww_1\biggl(\delta,\frac{d_1}{2}\biggr)
\le2 \rsww_1\biggl(\delta,\frac{\rho}{2}\biggr) ,
\]
which tends to zero with $\delta$ by Assumption \ref{eqwrsw},
and we finish the proof in case \hyperlink{caseii}{(2)}.

The proof for case \hyperlink{caseiii}{(3)} is easily obtained by considering the
dual closed crossing
from ${\partial}_1$ to ${\partial}_3$ and applying case
\hyperlink{caseii}{(2)}. The details are left to the reader.
\end{pf*}

The following lemma shows that it is unlikely to see three crossings
approaching the same boundary point.
This would make the crossing event unstable under boundary perturbation,
so in principle the lemma can be deduced from the previous one,
but instead we give a self-contained proof.
The lemma would follow from considering
the three arm event in \textit{half-annuli},
for whose probability Aizenman proposed an
argument to be comparable to~$(r/R)^2$.
We will use a version of his reasoning along with RSW
techniques to show a weaker estimate, sufficient for our purposes.
\begin{lemma}\label{lthree}
Let $Q$ be a quad with smooth sides, two opposite being labeled $\beta
$ and $\beta'$,
and with the tiles on each of the two other sides~$\nu$ and~$\nu'$
possibly declared all open or all closed.
Then there is a function $\rsww_Q(\delta)$, tending to zero with
$\delta$
and such that the following estimate holds.

Denote by $S(\theta)$ the event
that \textup{there are three crossings
with some prescribed alternating order
(say closed, open, and closed)
inside $[Q]$
between the set $\theta\subset\beta'$ and $\beta$}.
Then for $|\eta|<\delta$, we have
\[
\Pdisc_\eta\{\exists\theta\subset\beta' \dvtx S(\theta), \diam
(\theta)<\delta\} \le\rsww_Q(\delta) .
\]
%
\end{lemma}
\begin{pf}
Fix $\epsilon>0$.
We shall prove that for sufficiently small $\delta$
the probability of the event in question is less than $\epsilon$.

Since crossings can in principle use tiles on the sides $\nu$ and $\nu'$
(which are declared open or closed),
we have to treat $\theta$'s near the corners of $Q$ separately.
Let $z_1$ and $z_2$ be the endpoints of $\beta'$.
Since tiles on each of the arcs $\nu$ and $\nu'$ can be used by at
most one of the three crossings,
event $S(\theta)$ implies existence of at least one crossing
from $\theta$ to $\beta$ inside $[Q]$, not touching the sides.
By assumption (\ref{eqwrsw}), there is $r\in(0,\delta)$ such
that\vadjust{\goodbreak}
probability to have a crossing between $\bigcup_iB(z_i,2r)$ and $\beta$
is smaller than $\epsilon/2$:
\[
\Pdisc_\eta\biggl\{\exists\theta\subset\bigcup_iB(z_i,2r) \dvtx S(\theta)\biggr\}
\le\frac\epsilon2 .
\]
Set $\beta'':=\beta'\setminus\bigcup_iB(z_i,r)$, then
to prove lemma it remains to show that
%
%
\begin{equation}\label{eqsbetapp}
\Pdisc_\eta\{\exists\theta\subset\beta'' \dvtx S(\theta), \diam
(\theta)<\delta\} \le\frac\epsilon2 .
\end{equation}
Using assumption (\ref{eqwrsw}), choose $\rho\in(0,r/4)$ such that
$\rsww_1(2\rho,r/2)<1/2$.

Cover $\beta''$ by a finite number $n=n(\rho,\beta'')$ of
overlapping arcs $\theta_j$
of diameter at most $\rho$.
Let $\delta<\rho$ be so small that
$\rsww_1(4\delta,r)<\epsilon/(12n)$.
Since $\beta'$ is smooth, and decreasing $\delta$ if necessary,
we can cover each arc $\theta_j$ by arcs $\theta^l_j$ of length~$3\delta$, overlapping at most thrice,
and such that any set $\theta\subset\beta''$ of diameter less than
$\delta$
is entirely contained within one of the arcs $\theta_j^l$.\vspace*{1pt}

To deduce (\ref{eqsbetapp}) it is sufficient to show that, for a fixed $j$,
%
%
\begin{equation}\label{eqsgammaj}
\sum_l \Pdisc_\eta[S(\theta_j^l)] \le\frac{\epsilon
}{2n} .
\end{equation}
Fix $j$ and let $B(z,\rho)$ be a radius $\rho$ ball,
containing $\theta_j$.
The point $z_1$ splits ${\partial}Q\setminus\theta_j$ into two arcs:
$\lambda'\subset\beta'$ and a $\lambda\supset\beta$.
Both start at endpoints of $\theta_j$
inside $B(z,\rho)$ and end at least $r$-away.

Observe that event $S(\theta_j^l)$ is contained in the event
$S'(\theta_j^l)$ that
\textit{there are three alternating crossings inside $[Q]$
between $\theta$ and $\beta\cup\nu\cup\nu'$}.\vspace*{1pt}

Condition on the event $S'(\theta_j^l)$,
and, starting from $z_2$, denote the three such crossings,
closest to $z_2$,
by $\gamma_1$, $\gamma_2$, $\gamma_3$.
Namely, let $\gamma_1$ be the closed crossing between
$\theta$ and $\beta\cup\nu\cup\nu'$ inside $[Q]$, closest to $z_2$.
Let $M_1$ be the component of connectivity of $[Q]\setminus\gamma_1$,
not containing $z_2$.
Then $\gamma_1$ depends only on the percolation configuration in
$[Q]\setminus M_1$.

Now, take $\gamma_2$ be the open crossing between
$\theta$ and $\beta\cup\nu\cup\nu'$ inside $M_1$, closest to~$\gamma_1$.
Let $M_2$ be the component of connectivity of $M_2\setminus\gamma_2$,
not containing $\gamma_1$.
Then $\gamma_1$ and $\gamma_2$ depend only on the percolation
configuration in
$[Q]\setminus M_2$.

We define $\gamma_3$ and $M_3$ similarly,
and observe that $\gamma_1$, $\gamma_2$, $\gamma_3$
depend only on the percolation configuration in
$[Q]\setminus M_3$.
Therefore, given $M_3$ and the restriction
$\omega\llcorner([Q]\setminus M_3)$
(i.e., the restriction of percolation configuration
$\omega$ to $[Q]\setminus M_3$), 
the conditional law of the restriction
of $\omega$ to $M_3$ is unbiased.

Let $S^*=S^*(M_3)$ be the event
that \textit{there is a closed crossing $\gamma_4$
between~$\lambda$ and $\gamma_3$ inside $M_3$}.
Using that the restriction of $\omega$ to $M_3$
is unbiased, we can then write
%
%
\begin{equation}\label{eqsstar12}
\Pdisc_\eta\bigl[S^* \mid M_3, \omega\llcorner{[Q]\setminus M_3}\bigr] =
\Pdisc_\eta[S^*] > 1-\prsw^{1}_{\eta}(z,\rho,r) > \tfrac12 .
\end{equation}
Above we use that (in an unbiased percolation configuration),
if there is no open crossing of the annulus $A(z,\rho,r)$,
then by duality there is a closed circuit (a path going around the annulus),
whose intersection with $M_3$ would then contain a crossing required
for $S^*$.

Note that when $S'(\theta_j^l)$ and $S^*$ occur,
the crossing $\gamma_3$ from [$S(\theta_j^l)$]
and the crossing $\gamma_4$ (from $S^*$)
together give\vspace*{1pt} a closed crossing between $\theta$ and $\lambda$,
and so the following event $S''(\theta_j^l)$ occurs:
\textit{there is a closed crossing $\gamma_3'$ from $\theta$ to~$\lambda$},
\textit{an open crossing from $\gamma_2'$ from $\theta$ to $\beta\cup\nu
\cup\nu'$},
\textit{and a closed crossing from~$\gamma_1'$ from $\theta$ to $\beta\cup\nu\cup\nu'$}.
Then the estimate (\ref{eqsstar12}) can be rephrased as
\[
\Pdisc_\eta[S''(\theta_j^l) \mid S'(\theta_j^l)] \ge
\tfrac12 ,
\]
and therefore
%
%
\begin{equation}\label{eq2sm}
\Pdisc_\eta[S(\theta_j^l)] \le\Pdisc_\eta[S'(\theta_j^l)] \le2
\Pdisc_\eta[S''(\theta_j^l)] .
\end{equation}

Consider some percolation configuration $\omega$ in $[Q]$,
and assume that $S''(\theta_j^l)$ holds for some $\theta
_j^l\subset\theta_j$.

Choose inside $Q$ a closed crossing $\gamma_3'(\omega)$
between $\lambda'$ and $\theta_j$;
and an open crossing $\gamma_2'(\omega)$ between $\lambda$
which are the closest to each other
(and the point~$z_1$ separating $\lambda'$ from $\lambda$).
Denote by $L=L(\omega)$ the union of tiles in $\gamma_2'$, $\gamma_3'$,
and in the part of $[Q]$ between them.

Then $\gamma_2'$ and $\gamma_3'$ depend
only on the restriction of $\omega$ to $L$,
and the rest is unbiased.
Moreover, if $\gamma_2'$ ends at a tile on $\theta_j$,
then $\gamma_3'$ ends at a neighboring tile on $\theta_j$---otherwise
using duality
reasoning we can choose two closer crossings.
Denote by $x$ the only common vertex of these two tiles
lying on the boundary of $N$.
Note that 
each configuration $\omega$ has a unique such point $x=x(\omega)$
and it depends only on the restriction of $\omega$ to $N$.

Now for $S(\theta_j^l)$ to occur, besides $\gamma_2'$ and
$\gamma_3'$
we must also have a closed crossing $\gamma_1'$ from $\theta$ to
$\beta\cup\nu\cup\nu'$.
The crossing $\gamma_1'$ necessarily lies in $[Q]\setminus N$, which
is unbiased on $\omega\llcorner N$.
Note also that $\gamma_1'$ contains a crossing of the annulus
$A(x,2\delta,r)$,
whose center $x$ is a random point depending on $\omega\llcorner L$ only.
Summing it up, we can estimate
\[
\Pdisc_\eta\{\omega\dvtx S''(\theta_j^l) \mid x(\omega)=x\} \le
\prsw^1_\eta(x,4\delta,r)\le\rsww_1(4\delta,r)<\frac
{\epsilon}{12n} .
\]
Since every point $x$ is covered by at most three arcs $\theta_j^l$,
we recall (\ref{eq2sm})
and conclude that
\begin{eqnarray*}
\sum_l \Pdisc_\eta[S(\theta_j^l)]&\le&
\sum_l 2\Pdisc_\eta[S''(\theta_j^l)]\\
&=&\sum_{x}{\sum_l 2\Pdisc_\eta\{\omega\dvtx S(\theta_j^l) \mid x(\omega
)=x\}} \Pdisc_\eta\{\omega\dvtx x(\omega)=x\}\\
&\le&\sum_{x}\biggl(2\cdot3\cdot\frac{\epsilon}{12n}\biggr) \Pdisc_\eta
\{\omega\dvtx x(\omega)=x\}\\
&\le&\frac{\epsilon}{2n} \sum_{x} \Pdisc_\eta\{\omega\dvtx x(\omega
)=x\}
\le\frac{\epsilon}{2n} ,
\end{eqnarray*}
so estimate (\ref{eqsgammaj}) and the lemma follow.\vadjust{\goodbreak}
\end{pf}

\section{Multi-scale bound on the four-arm event on $\Z^2$ (by
Christophe Garban)}\label{appgarban}

In this Appendix, we give a proof of assumption (\ref{eqpivo}) for
critical bond percolation
on the square lattice. Namely:
%
\begin{lemma}\label{lemgarban}
For critical bond percolation
on the square lattice with mesh~$|\eta|$,
there is a positive $\epsilon$ such that
the probability of a four arm event satisfies
%
%
\begin{equation}
\label{eqpivosquare}
\prsw^{4}_{\eta}(z,r,R) \le\const\biggl(\frac{r}{R}\biggr)^{1+\epsilon} ,
\end{equation}
whenever $|\eta|<r$.
\end{lemma}

The case $r=|\eta|$ can be extracted from \cite{kesten-scaling} as
well as \cite{benjamini-kalai-schramm}
or \cite{schramm-steif-aom}. Note that the above multi-scale
generalization is not a consequence of the so-called
``quasi-multiplicativity'' property.
If one wanted to use this property, the problem would more or less boil
down to the \textit{existence} of the \textit{four-arm} critical exponent.
Since its existence is still open, we believe that Lemma~\ref{lemgarban} cannot be extracted directly from the above mentioned papers.

In the papers \cite{benjamini-kalai-schramm,schramm-steif-aom}, the
general idea which enables one to prove that points are unlikely
to be pivotal is the observation that a macroscopic crossing event for
critical percolation is asymptotically uncorrelated with the Majority
Boolean function defined on the same set of bits. This asymptotic
uncorrelation can be understood using an exploration path which will
compute the percolation event while
giving little information on the Majority event. A careful analysis
then shows that this decorrelation is possible only if points are
unlikely to be pivotal for the
percolation event. This heuristical program has been carried out in a
very convenient manner in \cite{odonnell-servedio}.

The present proof has the same flavor except that one now looks at
the correlation of a percolation event in a domain of diameter $O(R)$
with a \textit{two-layers majority function}:\vspace*{1pt} namely the bits of this
Majority function are now indexed by the $O(R^2/r^2)$ $r$-squares
included in the domain and for each such $r$-square, its
corresponding bit $\pm1$ depends on how ``connected'' the
percolation configuration is within this square. This setup allows to
``interpolate'' the above program from the macroscopic scale $R$ to
the mesoscopic scale~$r$.\looseness=1
\begin{pf*}{Proof of Lemma \ref{lemgarban}}
We leave out a few details, which are easy to fill in for readers
familiar with
the applications of the RSW theory, similar to
\cite{kesten-scaling,smirnov-werner,garban-pete-schramm}.

After rescaling and, if needed, changing the radii by bounded factors,
we can assume
without loss of generality that the lattice mesh is $1$ and $r$, $R$
are positive integers.


Denote by $Q$ an $R\times R$ square, and
cut the concentric $\frac13Q$ square
into $r\times r$ squares, denoted by $Q_j$, with $j=1,\ldots,(\frac{R}{3r})^2$.
Denote\vspace*{1pt} by $X=2\cdot1_{\cross{Q}}-1$ the ``crossing random variable''
equal to $1$ when
$Q$ is crossed and $-1$ otherwise.
Given a percolation configuration $\omega$, we say that
$Q_j$ is pivotal for $X$, if altering $\omega$
so that all the bonds in $Q_j$ are open,
and so that all the bonds in $Q_j$ are closed
yields two different values of $X$.
Note that $Q_j$ is pivotal for~$X$ if and only if
then there are two open arms connecting $\partial Q_j$ to
$\partial_0Q$ and $\partial_2Q$, and
two dual closed arms connecting $\partial Q_j$ to
$\partial_1Q$ and $\partial_3Q$.
By arm separation properties similar to ones used in \cite{smirnov-werner},
%
%
\begin{equation}\label{eqpivo4}
\prsw^{4}(r,R) \asymp\mathbf{P}[Q_j\mbox{ pivotal for }X] ,
\end{equation}
and so it is sufficient to estimate the latter probability,
or its sum over all~$Q_j$'s.

Let $S_j$ be an annulus $(1-\delta)Q_j \setminus(1-2\delta)Q_j$
where $ \alpha Q_j$ denotes a copy shrank of $Q_j$ by a factor of
$\alpha$ (if $r$ is not exactly divisible by $\alpha$, the copy
shrank is understood modulo integer parts).
The fixed parameter $\delta>0$ will be chosen later. Let $S_j^*$
denote the same annulus but on the dual lattice: that is, $S_j^* = S_j
+ (1/2,1/2)$.


Define a random variable $C_j$ equal to:
\begin{itemize}
\item$1$, if there is an open circuit in $S_j$ and no dual-closed
circuit in $S_j^*$,
\item$-1$, if there is a dual-closed circuit in $S_j^*$ and no open
circuit in $S_j$,
\item$0$, otherwise.\vspace*{5pt}
\end{itemize}
Note that by symmetry and by the RSW theory, up to the constant
depending on $\delta$,
%
%
\begin{equation}\label{eqcj}
\E[C_j] =0 ,\qquad \E[C_j^2] \asymp1 .
\end{equation}

Now, let us argue that if $\delta>0$ is chosen small enough, one has
%
%
\begin{equation}\label{eqsepar}
\E[X C_j] \asymp\mathbf{P}[Q_j\mbox{ pivotal for }X] \qquad[ \mbox{$\asymp$}\prsw
^{4}(r,R) ].
\end{equation}

Indeed, one has
%
%
\begin{equation}\label{eqXCj}
\E[X C_j] = \mathbf{P}[Q_j\mbox{ pivotal for }X] \E[X C_j \mid
Q_j\mbox{ pivotal for }X] ,
\end{equation}
since, conditioned on the event that $Q_j$ is not pivotal for $X$,
$C_j$ is independent of $X$ and is such that its (conditional)
expectation is still 0. Therefore, it remains to bound from below $\E
[X C_j \mid Q_j\mbox{ pivotal for }X]$. Let $\rho=\rho(\delta
)>0$ be such that
$\mathbf{P}[C_j =1] = \mathbf{P}[C_j=-1] > \rho$.
\begin{eqnarray*}
&&\E[X C_j \mid Q_j\mbox{ pivotal for }X] \\
&&\qquad = \mathbf{P}[C_j =1] \E
[X \mid C_j=1 \mbox{ and } Q_j\mbox{ pivotal for }X] \\
&&\qquad\quad{} - \mathbf{P}[C_j =1] \E[X \mid C_j=-1 \mbox{ and } Q_j\mbox{
pivotal for }X] .
\end{eqnarray*}

Now, again by arm separation properties similar to ones used in
\cite{smirnov-werner,schramm-steif-aom},
and using the important fact that the event $\{ C_j=1 \}$ is increasing
which enables to use FKG for the given conditional law inside $Q_j$,
one can see that,
if $\delta$ is chosen small enough,
\[
\E[X \mid C_j=1 \mbox{ and } Q_j\mbox{ pivotal for }X] > \tfrac1
4 .
\]
One has the opposite bound for the term conditioned on $\{ C_j =-1 \}
$. All together this gives
\[
\E[X C_j \mid Q_j\mbox{ pivotal for }X] > \frac{\rho} 2 ,
\]
which implies the desired estimate (\ref{eqsepar}) if $\delta$ is
chosen to be small enough.

Consider $Q$ with open boundary conditions on the side $\partial_0 Q$
and dual closed on the complementary three sides.
Let $\interface$ be the interface running between the two ends
of the side $\partial_0Q$ and separating the open cluster
rooted on it from the dual closed cluster rooted on the three other sides.

Denote by $Y_j$ the event that
$\interface$ intersects $Q_j$,
as well as its indicator function.
Note that, by the RSW theory,
%
%
\begin{equation}\label{eqinter}
\E[Y_j] \lesssim\biggl(\frac rR\biggr)^{2\epsilon}
\end{equation}
for some positive $\epsilon$.
Note that the interface $\interface$ drawn from $Q(0,0)$ until some
stopping time
depends only on the percolation configuration in the
immediate neighborhood of the drawn part.
Drawing $\interface$ until the first time it hits~$\partial Q_j$,
we conclude that $Y_j$ is independent from the inside of $Q_j$ and
hence from $C_j$.

Note also that, if $Q_j$ is pivotal for $X$,
then there are four alternating arms connecting $\partial Q_j$ to
four sides of $Q$ which forces the interface $\interface$ to intersect~$Q_j$, and $Y_j$ to occur.
Therefore, similarly as for the estimate (\ref{eqXCj}) above, one can
check that
%
%
\begin{equation}
\E[XC_j] = \E[XC_jY_j] ,
\end{equation}
which combined with (\ref{eqsepar}) gives
%
%
\begin{equation}\label{eqpivoexp}
\mathbf{P}[Q_j\mbox{ pivotal for }X] \asymp\E[XC_jY_j] .
\end{equation}
Summing over all $Q_j$'s, we
use the Cauchy--Schwarz inequality to write
%
%
\begin{eqnarray}\label{eqdiag}
\sum_j\E[XC_jY_j] &=& \E\biggl[X\tsum_j{C_jY_j}\biggr] \le\sqrt{\E[X^2]\cdot
\E\biggl[\biggl(\tsum_j{C_jY_j}\biggr)^2\biggr]}\nonumber\\[-8pt]\\[-8pt]
&=& \sqrt{\E[1]\cdot\E\biggl[{\tsum_{i,j}{C_iY_iC_jY_j}}\biggr]} = \sqrt
{{\tsum_{i,j}\E[{C_iY_iC_jY_j}]}}.\nonumber
\end{eqnarray}
In the last sum, nondiagonal terms vanish.
Indeed, let $i\neq j$ and stop the~in\-terface $\interface$ the first
time it touched both squares
$Q_i$ and $Q_j$, or when it~ends.

Denote by $Q_k$ the first square of $Q_i$, $Q_j$ to be hit
(or $Q_i$ if none was)
and by $Q_l$ the other one,
so that $(k,l)$ is a superposition of $(i,j)$.
Let $W$ be the percolation configuration in the immediate neighborhood
of the interface so far, and in $Q_k$.\vadjust{\eject}

Then $W$ determines $Y_i$, $Y_j$, $C_k$, while being
independent of $C_l$. Thus
by the identity in (\ref{eqcj}),
\[
\E[C_iY_iC_jY_j\mid W] 
= \E[Y_iY_jC_k\mid W]\cdot\E[C_l] = 0.
\]
By the total expectation formula we conclude that,
for $i\neq j$,
\[
\E[C_iY_iC_jY_j] = 0 .
\]
Thus, we can continue (\ref{eqdiag}), leaving only the diagonal terms,
and use (\ref{eqinter}) and (\ref{eqcj}) to write
\begin{eqnarray*}
\cdots&=& \sqrt{{\tsum_{j}{\E[C_j^2Y_j^2]}}} = \sqrt{{\tsum
_{j}{\E[C_j^2]\cdot\E[Y_j^2]}}}\\
&\lesssim&\sqrt{\tsum_{j}1\cdot\biggl(\frac rR\biggr)^{2\epsilon}}
= \biggl(\frac rR\biggr)^{\epsilon-1} .
\end{eqnarray*}
Combining this with (\ref{eqpivo4}), (\ref{eqpivoexp}), (\ref{eqdiag}),
we conclude
\begin{eqnarray*}
\prsw^{4}(r,R)&\lesssim&
\biggl(\frac{3r}{R}\biggr)^2 \tsum_j\mathbf{P}[Q_j\mbox{ pivotal for }X]\\
&\lesssim&
\biggl(\frac{r}{R}\biggr)^2\biggl(\frac rR\biggr)^{\epsilon-1} = \biggl(\frac rR\biggr)^{1+\epsilon} ,
\end{eqnarray*}
proving the lemma.
\end{pf*}
\end{appendix}

\section*{Acknowledgments}
We are thankful to Christophe Garban and Wendelin Werner for useful
conversations. The second author is grateful to Hugo Duminil-Copin,
Cl\'{e}ment Hongler, G\'{a}bor Pete, Jeff Steif, the referee, and
especially Christophe Garban and Boris Tsirelson for their comments on
the manuscript.


%
%

%
\printaddresses

\end{document}